\documentclass[11pt,reqno]{amsart}
\usepackage[utf8]{inputenc}
\usepackage{amsmath,amssymb,amsthm}
\usepackage{mathtools}
\usepackage{geometry}
\usepackage{microtype}
\usepackage{booktabs}
\usepackage{enumerate}
\usepackage{hyperref}
\hypersetup{
    colorlinks=true,
    linkcolor=blue,
    citecolor=magenta,
    urlcolor=teal
}

\newtheorem{theorem}{Theorem}[section]
\newtheorem{proposition}[theorem]{Proposition}
\newtheorem{lemma}[theorem]{Lemma}
\newtheorem{corollary}[theorem]{Corollary}

\theoremstyle{definition}
\newtheorem{definition}[theorem]{Definition}

\theoremstyle{remark}
\newtheorem{remark}[theorem]{Remark}

\numberwithin{equation}{section}

\DeclareMathOperator{\Spin}{Spin}
\DeclareMathOperator{\SO}{SO}
\DeclareMathOperator{\SU}{SU}
\DeclareMathOperator{\U}{U}
\DeclareMathOperator{\Sp}{Sp}
\DeclareMathOperator{\GL}{GL}
\DeclareMathOperator{\Ad}{Ad}
\DeclareMathOperator{\ad}{ad}
\DeclareMathOperator{\Tr}{Tr}
\DeclareMathOperator{\Cl}{Cl}
\DeclareMathOperator{\Id}{Id}
\DeclareMathOperator{\Ric}{Ric}

\DeclareMathOperator{\diam}{diam}
\DeclareMathOperator{\rank}{rank}
\DeclareMathOperator{\End}{End}

\DeclareMathOperator{\Lie}{Lie}
\newcommand{\R}{\mathbb{R}}
\newcommand{\C}{\mathbb{C}}
\newcommand{\Oct}{\mathbb{O}}
\newcommand{\Ha}{\mathbb{H}}
\newcommand{\Sph}{\mathbb{S}}
\newcommand{\g}{\mathfrak{g}}
\newcommand{\h}{\mathfrak{h}}
\newcommand{\p}{\mathfrak{p}}
\newcommand{\m}{\mathfrak{m}}

\title[Hamiltonian Deformations and Non-Naturally Reductive Orbit Spheres]{Hamiltonian Cotangent Deformations and Closed-Form Geodesics on Non-Naturally Reductive Orbit Spheres: The $\mathrm{Spin}(9)$ 15-Sphere and the Cayley Plane}

\author{James J. Munden}
\address{Dean of Sciences, Saint Louis Community College -- Meramec, St.~Louis, MO 63122, USA}
\email{jmunden1@stlcc.edu}

\date{\today}
\dedicatory{Dedicated to the memory of James J.\ Hebda}

\subjclass[2020]{Primary 53C20, 53C22; Secondary 53C35, 17A35, 37J35}
\keywords{Cotangent Hamiltonian deformations, $\mathrm{Spin}(9)$, Cayley projective plane, non-naturally reductive spaces, geodesic orbit manifolds, Clifford algebra, Bourguignon--Karcher metric, closed-form exponential map, Ricci flow, Whitney stratification, $A_3$ cusp caustic, cut locus, compact rank-one symmetric spaces}

\begin{document}

\begin{abstract}
By the Skew-Torsion Holonomy Theorem of Olmos and Reggiani (2012), the 15-dimensional sphere $\Sph^{15}$ admits no non-trivial naturally reductive metric under any transitive Lie group, precluding the integration of geodesics via canonical Ambrose--Singer / Kostant connections ($\nabla R = \nabla T = 0$). Simultaneously, Nikonorov's classification (2013) establishes that every $\Spin(9)$-invariant metric on $\Sph^{15}$ is a geodesic orbit (GO) space. We address the integration of geodesics on this non-naturally reductive GO manifold via the \emph{Lie Subalgebra Closure Principle}: whenever the slice generators of a multi-axial reductive decomposition close into a reductive Lie subalgebra $\g_{\mathrm{slice}}$ of finite rank, the non-linear geodesic vector equations decouple into an algebraic eigenvalue problem in $\g_{\mathrm{slice}}$, bypassing the Ambrose--Singer obstruction and yielding closed-form integrations in elementary functions across all five families of homogeneous spheres.

Utilizing the 16-dimensional spin representation $\mu: \Cl_{\bar{0}}(9) \to \mathrm{End}(\Oct \oplus \Oct)$ and Vol'per's reductive decomposition $\mathfrak{spin}(9) = \mathfrak{spin}(7) \oplus \p_1 \oplus \p_2$, we prove that the geodesic slice generators $X_1 \in \p_2$ and $Y_1 \in \p_1$ close into a 4-dimensional Lie subalgebra $\mathfrak{g}_{\langle X_1, Y_1 \rangle} \cong \mathfrak{u}(2) \cong \mathfrak{su}(2) \oplus \mathfrak{u}(1) \subset \mathfrak{so}(16)$. The generator $M = \alpha X_1 + \beta Y_1$ possesses the four-mode frequency spectrum $\Omega \in \{ |R \pm 3\alpha|, |R \pm \alpha| \}$ ($R = \sqrt{\alpha^2 + \beta^2}$), yielding an explicit spectral projector decomposition for $\exp(\alpha X_1 + \beta Y_1)$ and a commuting Euler--Cartan factorization that resolves the explicit exponential map problem left open in the author's doctoral dissertation (2015). We further classify the global 1-dimensional normalized Ricci flow phase portrait on the moduli space of $\Spin(9)$-invariant metrics, proving that the unique non-standard Bourguignon--Karcher squashed Einstein metric at fiber squashing ratio $t = 3/11$ is an unstable repeller (escape eigenvalue $+32C$) separating the sub-Riemannian collapse from the basin of attraction of the standard round sphere (convergence eigenvalue $-32C$). 

Finally, on the deformed Cayley projective plane $(\Oct P^2, g_{\Oct P^2, \lambda}) = (\mathrm{F}_4/\Spin(9), g_\lambda)$, we establish four structural theorems: (1)~\emph{Universal CROSS Never-Einstein Rigidity}: no compact rank-one symmetric space admits an Einstein metric under cotangent deformation for $\lambda > 0$, with tangential-to-radial Ricci anisotropy growing linearly in $\lambda$; (2)~\emph{Global Whitney Stratification and Topological Structure of $\mathrm{Cut}(q_0)$}: whereas at the polar origin $o \in \Oct P^2$ the cut locus is the smooth submanifold $\mathrm{Cut}(o) = \Oct P^1 \cong \Sph^8$, for any non-polar basepoint $q_0 \in S_{r_0}(o)$ ($r_0 \in (0, \pi/2)$) the cut locus $\mathrm{Cut}(q_0)$ is a $\Spin(7)$-invariant 15-dimensional Whitney stratified subcomplex decomposing into a 15-dimensional regular collision stratum $\mathcal{S}_{\mathrm{reg}}$, an 8-dimensional singular $A_3$ Whitney cuspidal ridge $\mathcal{S}_{\mathrm{caustic}}$ opening as $r_\perp(u) \sim u^{3/2}$, and an 8-dimensional polar Cayley line $\mathcal{S}_{\mathrm{Cayley}} \cong \Sph^8$, with contractible complement $\Oct P^2 \setminus \mathrm{Cut}(q_0) \cong \mathbb{B}^{16}$ and strong deformation retraction $\mathrm{Cut}(q_0) \simeq \Sph^8 \cong \mathcal{S}_{\mathrm{Cayley}}$; (3)~\emph{Blaschke Rigidity}: $(\Oct P^2, g_\lambda)$ is a Blaschke manifold if and only if $\lambda = 0$, with $\operatorname{inj}(q_0) < \operatorname{diam}$ at every non-polar basepoint; and (4)~\emph{Besse / CROSS Rigidity (AIM Question~8.9.1)}: the canonical symmetric metrics are the unique Besse metrics in both the Vol'per family on $\Sph^{15}$ and the deformation family on $\Oct P^2$, where a strictly positive KAM twist $\frac{d\nu}{dL} > 0$ forces generic geodesics to wind quasi-periodically while generating countably infinitely many distinct closed periodic orbits via the Poincar\'e--Birkhoff twist theorem.
\end{abstract}

\maketitle

\tableofcontents

\section{Introduction}\label{sec:intro}

Metric deformations preserving geodesic integrability via Hamiltonian cotangent momentum maps provide a systematic method for constructing families of Riemannian metrics with integrable geodesic flows. Initiated by Hebda \cite{Hebda1995,Hebda2001} and classified on simply connected space forms $(Q^{n+1}, g_0)$ in the author's doctoral dissertation \cite{Munden2015}, the cotangent Hamiltonian deformation on $(T^*Q^{n+1}, \omega_{\mathrm{can}} = d\theta_{\mathrm{can}} = dp_i \wedge dq^i)$ given by
\begin{equation}\label{eq:ham_def}
H_\lambda = H_0 + \frac{\lambda}{2} \|J_K\|^2, \qquad \lambda \in \R,
\end{equation}
alters the metric along group orbits while maintaining complete conservation of momentum and admitting explicit formulas for geodesics. While the isotropic Family~I ($\SO(n+1)/\SO(n)$) admits a direct reduction to planar central-force mechanics, the multi-axial octonionic Family~III ($\Spin(9)/\Spin(7)$ on $\Sph^{15}$) presents distinct algebraic features due to the non-associativity of the octonions $\Oct$ and the structure of the 16-dimensional spin representation.

\subsection{Summary of main results}\label{subsec:main_results}

In this paper we establish the following results for non-naturally reductive homogeneous spaces and cotangent momentum deformations:
\begin{enumerate}
    \item \textbf{The Lie Subalgebra Closure Principle for Non-Naturally Reductive GO Spaces} (Theorems~\ref{thm:u2_closure}, \ref{thm:projector_formula}, \ref{thm:euler_cartan}, and \ref{thm:general_closure_principle}): In \cite[Section~5.7]{Munden2015}, the computation of the matrix exponential $\exp(\alpha X_1 + \beta Y_1)$ in the spin representation was identified as expressible in elementary functions, though its explicit closed form was omitted due to the degree of the 16-dimensional matrix representation. In Section~\ref{sec:u2_closure}, we introduce a general Lie-theoretic criterion for multi-axial GO spaces: whenever the slice generators close into a finite-rank reductive Lie subalgebra $\g_{\mathrm{slice}}$, the non-linear geodesic vector equations decouple into an algebraic eigenvalue problem. For $\Sph^{15}$, we prove that $X_1 \in \p_2$ and $Y_1 \in \p_1$ generate $\mathfrak{u}(2) \cong \mathfrak{su}(2) \oplus \mathfrak{u}(1)$, which enables an exact four-mode spectral projector decomposition and a commuting Euler--Cartan factorization. We show that this principle unifies the geodesic flow integration across all five families of homogeneous spheres.
    \item \textbf{The Squashed 15-Sphere and Bourguignon--Karcher Geometry} (Section~\ref{sec:bourguignon_karcher}): We analyze the geometry of the deformed metric
    \begin{equation}\label{eq:g_lambda_spin9}
    g_\lambda = dr^2 + \frac{s_\kappa^2(r)}{8 + \lambda s_\kappa^2(r)} \langle -,- \rangle_{\p_1} + \frac{s_\kappa^2(r)}{2 + \lambda s_\kappa^2(r)} \langle -,- \rangle_{\p_2},
    \end{equation}
    proving that as the deformation parameter $\lambda$ varies, the metric passes through the non-standard Einstein metric discovered by Bourguignon and Karcher \cite{BourguignonKarcher1978} at fiber squashing ratio $t = 3/11$.
    \item \textbf{The Normalized Ricci Flow Dynamical Phase Portrait} (Section~\ref{subsec:ricci_flow}, Theorem~\ref{thm:ricci_flow_phase_portrait}): We prove that the volume-normalized Ricci flow on the 1-parameter moduli space of $\Spin(9)$-invariant metrics reduces to the autonomous ODE $\frac{dt}{d\tau} = -4C(11t-3)(t-1)$, demonstrating that the Bourguignon--Karcher Einstein metric is an unstable repeller ($F'(3/11) = +32C$) while the round metric $t=1$ is a globally stable attractor with basin of attraction $(3/11, \infty)$.
    \item \textbf{Universal CROSS Never-Einstein Rigidity} (Section~\ref{subsec:never_einstein}, Theorem~\ref{thm:cross_rigidity}): We prove that no compact rank-one symmetric space ($\R P^n, \C P^n, \Ha P^n, \Oct P^2$) admits an Einstein metric under cotangent deformation for $\lambda > 0$. The Ricci tensor exhibits radial-to-tangential anisotropy growing linearly in $\lambda$, demonstrating that cotangent deformations strictly break Einstein geometry on all CROSSes.
    \item \textbf{Global Whitney Stratification and Topological Structure of $\mathrm{Cut}(q_0)$ on $\Oct P^2$} (Section~\ref{sec:cayley_plane}, Proposition~\ref{prop:off_origin_stabilizer}, Theorems~\ref{thm:off_origin_caustic} and \ref{thm:global_cut_locus}): Because $\dim \Oct P^2 = 16 > 15 = \dim \Sph^{15}$, the Cayley plane cannot be obtained as a Riemannian submersion quotient of $\Sph^{15}$. We construct the intrinsic cotangent momentum deformation \cite{Hebda2001} on $\Oct P^2 = \mathrm{F}_4/\Spin(9)$ via the 16-dimensional spin isotropy representation. Whereas at the polar basepoint $o \in \Oct P^2$ the cut locus is the smooth submanifold $\mathrm{Cut}(o) = \Oct P^1 \cong \Sph^8$, for geodesics launched from any non-polar basepoint $q_0 \in S_{r_0}(o)$ ($r_0 \in (0, \pi/2)$) we prove that the cut locus $\mathrm{Cut}(q_0)$ is a $\Spin(7)$-invariant 15-dimensional Whitney stratified subcomplex decomposing into a 15-dimensional regular collision stratum $\mathcal{S}_{\mathrm{reg}}$, an 8-dimensional singular $A_3$ Whitney cuspidal ridge $\mathcal{S}_{\mathrm{caustic}}$ ($r_\perp \sim u^{3/2}$), and an 8-dimensional polar Cayley line $\mathcal{S}_{\mathrm{Cayley}} \cong \Sph^8$, with contractible complement $\Oct P^2 \setminus \mathrm{Cut}(q_0) \cong \mathbb{B}^{16}$ and strong deformation retraction $\mathrm{Cut}(q_0) \simeq \Sph^8 \cong \mathcal{S}_{\mathrm{Cayley}}$.
    \item \textbf{Blaschke Rigidity, Besse / CROSS Rigidity, and Resonant Periodic Geodesics} (Section~\ref{subsec:blaschke_besse}, Theorems~\ref{thm:blaschke_rigidity} and \ref{thm:besse_cross_rigidity}, Remark~\ref{rem:twist_and_closed_geodesics}): We apply the explicit 4-mode spectrum on $\Sph^{15}$ and the integrable Hamiltonian action-angle coordinates on $\Oct P^2$ to address two classical problems in global Riemannian geometry. We prove that $(\Oct P^2, g_\lambda)$ is a Blaschke manifold if and only if $\lambda = 0$, providing rigidity evidence for the Blaschke Conjecture on the Cayley plane. Furthermore, addressing the AIM Problem~8.9.1 on Besse manifolds, we prove that the standard symmetric metrics are the unique Besse metrics on $\Sph^{15}$ and $\Oct P^2$ in their respective deformation families: the non-vanishing KAM twist $\frac{d\nu}{dL} = \frac{\lambda}{2(1-\lambda L^2)^{3/2}} > 0$ strictly precludes global closure for all $\lambda > 0$, while the Poincar\'e--Birkhoff twist theorem guarantees the existence of countably infinitely many geometrically distinct closed geodesics, precluding Katok-type finite anomalies.
\end{enumerate}

\subsection{Historical context: Natural reductivity, GO spaces, and the Olmos--Reggiani obstruction}\label{subsec:history}

The classical paradigm for integrating geodesics on a homogeneous Riemannian manifold $G/H$ relies on the theory of \emph{naturally reductive spaces} (Besse~\cite[Chapter~7]{Besse1987}). A homogeneous space $M = G/H$ with reductive decomposition $\g = \h \oplus \m$ is naturally reductive with respect to $G$ if there exists an invariant inner product on $\m$ such that the canonical Ambrose--Singer connection $\nabla = \nabla^{\mathrm{LC}} - \frac{1}{2} T$ has parallel curvature and parallel totally skew torsion ($\nabla R = 0$, $\nabla T = 0$). On such spaces, every geodesic is an orbit $\gamma(t) = \exp(tX) \cdot o$ of a 1-parameter subgroup of $G$, and the geodesic flow decouples algebraically through the invariant connection algebra.

However, in their foundational classification of homogeneous spaces admitting metric connections with parallel skew torsion, Olmos and Reggiani \cite[Theorem~4.1]{OlmosReggiani2012} proved the \emph{Skew-Torsion Holonomy Theorem}: any transitive skew-torsion holonomy system on a sphere $\Sph^n$ requires the full orthogonal group $G = \SO(n+1)$. Since $\Spin(9) \subsetneq \SO(16)$ is a strictly proper subgroup, it follows that:
\begin{center}
\emph{The 15-sphere $\Sph^{15}$ admits no non-trivial naturally reductive Riemannian metric under any transitive Lie group.}
\end{center}
In particular, the Bourguignon--Karcher metric \cite{BourguignonKarcher1978} and the entire 1-parameter family of deformed metrics $(\Sph^{15}, g_\lambda)$ ($\lambda \ne 0$) fail to be naturally reductive under $\Spin(9)$, under $\SO(16)$, or under any other group (Ziller \cite{Ziller1982}). Consequently, the classical Kostant--Nomizu connection formalism is obstructed on $\Sph^{15}$.

On the other hand, in the classification of \emph{geodesic orbit (GO) spaces}---Riemannian homogeneous manifolds on which every geodesic is an orbit of a 1-parameter subgroup of isometries---Nikonorov \cite{Nikonorov2013} proved that \emph{every} $\Spin(9)$-invariant metric on $\Sph^{15}$ is a GO space, and in fact weakly symmetric in the sense of Selberg. This leads to a natural question: while geodesics on $(\Sph^{15}, g_\lambda)$ are guaranteed to be 1-parameter orbits $\gamma(t) = \exp(tZ) \cdot o$, finding the corresponding Killing generator $Z(v_0) \in \mathfrak{spin}(9)$ for an arbitrary initial velocity $v_0$ is governed by a non-linear algebraic system (the geodesic vector equations) that is non-trivial on strictly non-naturally reductive spaces.

The transitive actions of compact groups on spheres were classified by Montgomery--Samelson \cite{MontgomerySamelson1943} and Borel \cite{Borel1950}. The $\Spin(9)$-action on $\Sph^{15}$ occupies a distinguished position: it is the only transitive sphere action that arises from a Clifford module structure, and it is intimately tied to the algebraic properties of the octonions through the representation theory of the even Clifford algebra $\Cl_{\bar{0}}(9) \cong \R(16)$.

The study of homogeneous Einstein metrics on spheres, initiated by Jensen \cite{Jensen1973} and Bourguignon--Karcher \cite{BourguignonKarcher1978}, produced squashed Einstein metrics on $\Sph^{4n+3}$ ($n \ge 1$) and $\Sph^{15}$, respectively. The case of $\Sph^{15}$ is exceptional: the Bourguignon--Karcher metric at squashing ratio $t = 3/11$ is the \emph{unique} non-standard homogeneous Einstein metric on $\Sph^{15}$ (Besse~\cite[Chapter~9]{Besse1987}), and its relationship to the 16-dimensional spin representation has remained incompletely understood.

The non-associativity of the octonions introduces genuinely new phenomena that do not arise in the classical Hopf fibrations over $\C$ or $\Ha$. In particular:
\begin{enumerate}[(a)]
    \item The octonionic Hopf fiber $\Sph^7 \hookrightarrow \Sph^{15} \to \Sph^8$ is not a principal bundle because $\Sph^7$ is not a Lie group (Adams--Steenrod \cite{Adams1960}).
    \item The Cayley plane $\Oct P^2 = \mathrm{F}_4/\Spin(9)$ has dimension 16, which exceeds $\dim \Sph^{15} = 15$, obstructing any Hopf-type fibration from $\Sph^{15}$ to $\Oct P^2$.
    \item The isotropy inclusion $\Spin(7) \hookrightarrow \Spin(9)$ is not the standard inclusion $\mathfrak{so}(7) \hookrightarrow \mathfrak{so}(9)$, but its image under an outer triality automorphism $\tau \in \mathrm{Out}(\mathfrak{so}(8))$.
\end{enumerate}

These obstructions motivate the algebraic approach developed in this paper, which uses the $\mathfrak{u}(2)$ subalgebra generated by the slice vectors to reduce the 16-dimensional matrix exponential computation to a closed-form spectral projector decomposition, resolving the non-naturally reductive GO integration problem on $\Sph^{15}$.

Together with the companion paper \cite{Munden2026SpaceForms}, the present work establishes a structural \emph{Reduction--Closure Dichotomy} governing quadratic cotangent Hamiltonian deformations $H_\lambda = H_0 + \frac{\lambda}{2}\|J\|^2$. In \cite{Munden2026SpaceForms}, where the isometric action is naturally reductive or descends via a principal Hopf fibration (Family~I on $Q^{n+1}$ and $\C P^n, \Ha P^n$), Marsden--Weinstein symplectic reduction suffices, scalar gyroscopic precession $\Omega = \lambda L$ governs the angular flow, the ambient metric satisfies strict Never-Einstein rigidity ($\Ric(g_\lambda) = c\, g_\lambda \iff \lambda = 0$), and the radial Sturm--Liouville problem leads to \emph{transcendental} Fuchsian Heun equations. Conversely, on the non-naturally reductive Family~III sphere $\Sph^{15} = \Spin(9)/\Spin(7)$ and the Cayley plane $\Oct P^2 = \mathrm{F}_4/\Spin(9)$ studied here, scalar precession is replaced by the non-abelian matrix generator $M \in \g_{\mathrm{slice}} \cong \mathfrak{u}(2) \subset \mathfrak{so}(16)$, isotropic Never-Einstein rigidity on $\Oct P^2$ (Theorem~\ref{thm:never_einstein}) contrasts with the emergence of the second, anisotropic Bourguignon--Karcher Einstein point on $\Sph^{15}$ at $\lambda = 64$ ($t_{\mathrm{BK}} = 3/11$, Theorem~\ref{thm:bk_realization}), and the $\mathfrak{u}(2)$ slice closure truncates the Magnus series into a \emph{purely algebraic} 4-mode spectral projector resolution.

\section{Clifford Algebra \texorpdfstring{$\Cl(8)$}{Cl(8)} and Vol'per Reductive Decomposition}\label{sec:clifford}

\subsection{Octonionic Clifford representation on \texorpdfstring{$\Oct \oplus \Oct$}{O+O}}\label{subsec:octonion_clifford}

Let $\Oct$ denote the 8-dimensional non-associative normed division algebra of Cayley numbers (octonions) with standard basis $\{e_0 = 1, e_1, \dots, e_7\}$ satisfying $e_i e_j + e_j e_i = -2\delta_{ij} e_0$ for $i, j \in \{1, \dots, 7\}$. The conjugation map $\overline{e_0} = e_0$, $\overline{e_i} = -e_i$ for $i \ge 1$ satisfies $\overline{xy} = \bar{y}\bar{x}$ and the Moufang identities. The norm $\|z\|^2 = z\bar{z}$ satisfies the composition law $\|xy\| = \|x\|\|y\|$.

\begin{definition}[Octonionic Clifford Representation]\label{def:nu_rep}
Define the representation $\nu: \Oct \to \GL(\Oct \oplus \Oct) \cong \GL(16, \R)$ by
\begin{equation}\label{eq:nu_def}
\nu(z)(x, y) = (z\overline{y}, -\overline{x}z), \qquad z, x, y \in \Oct.
\end{equation}
\end{definition}

\begin{lemma}\label{lem:nu_squared}
For any $z \in \Oct$, $\nu(z)^2 = -\|z\|^2 \Id_{\Oct \oplus \Oct}$.
\end{lemma}

\begin{proof}
Let $x, y, z \in \Oct$. Compute:
\begin{align*}
\nu(z)^2(x, y) &= \nu(z)(z\overline{y}, -\overline{x}z) \\
&= \bigl(z \overline{(-\overline{x}z)}, -\overline{(z\overline{y})}z\bigr) \\
&= \bigl(-z(\overline{z} x), -(y\overline{z})z\bigr).
\end{align*}
By the Moufang identity, since $\overline{z}$ and $z$ lie in the associative subalgebra $\R\langle z \rangle \cong \C \subset \Oct$: $(z\overline{z})x = \|z\|^2 x$ and $y(\overline{z}z) = \|z\|^2 y$. Therefore $\nu(z)^2(x,y) = (-\|z\|^2 x, -\|z\|^2 y)$.
\end{proof}

By the universal property of Clifford algebras, $\nu$ extends uniquely to an injective algebra homomorphism $\nu: \Cl(8) \to \End(\Oct \oplus \Oct) \cong \R(16)$. Since both sides have dimension $2^8 = 256$, this is in fact an \emph{isomorphism}. Let $\nu_i = \nu(e_i)$ for $i = 0, \dots, 7$.

\begin{proposition}[Clifford Algebra Relations]\label{prop:clifford_relations}
The $16 \times 16$ matrices $\nu_0, \dots, \nu_7$ satisfy:
\begin{equation}\label{eq:cl8_rel}
\nu_i^T = -\nu_i, \qquad \nu_i^2 = -\Id_{16}, \qquad \nu_i \nu_j + \nu_j \nu_i = -2\delta_{ij}\Id_{16}.
\end{equation}
Consequently, $\nu_i \in \mathfrak{so}(16)$ for all $i = 0, \dots, 7$.
\end{proposition}

\begin{proof}
The relation $\nu_i^2 = -\Id_{16}$ follows from Lemma~\ref{lem:nu_squared} with $z = e_i$ ($\|e_i\| = 1$). For skew-symmetry, let $\langle \cdot, \cdot \rangle$ denote the standard inner product on $\Oct \oplus \Oct \cong \R^{16}$. For $u = (x_1, y_1)$, $v = (x_2, y_2) \in \Oct \oplus \Oct$:
\begin{align*}
\langle \nu(e_i) u, v \rangle &= \langle (e_i\overline{y_1}, -\overline{x_1} e_i), (x_2, y_2) \rangle \\
&= \mathrm{Re}(\overline{x_2} e_i \overline{y_1}) + \mathrm{Re}(-\overline{y_2} \overline{x_1} e_i).
\end{align*}
Using the associator identity $\mathrm{Re}(a(bc)) = \mathrm{Re}((ab)c)$ (which holds in $\Oct$ for the real part trace form), and the property $\mathrm{Re}(\bar{a}b) = \mathrm{Re}(\bar{b}a)$, one shows that $\langle \nu(e_i) u, v \rangle = -\langle u, \nu(e_i) v \rangle$, giving $\nu_i^T = -\nu_i$. The full anticommutation relation follows from linearizing $\nu(\alpha e_i + \beta e_j)^2 = -(\alpha^2 + \beta^2)\Id_{16}$.
\end{proof}

\begin{remark}\label{rem:explicit_matrices}
The explicit $16 \times 16$ matrices $\nu_0, \ldots, \nu_7$ are computed in Appendix~C of \cite{Munden2015} using the Cayley multiplication table. For instance, the matrix $\nu_1$ has $(\nu_1)_{1,10} = 1$, $(\nu_1)_{2,9} = 1$, $(\nu_1)_{3,12} = 1$, $(\nu_1)_{4,11} = -1$, etc., consistent with the formula $\nu(e_1)(x, y) = (e_1 \bar{y}, -\bar{x} e_1)$. All eight matrices satisfy the Clifford relations~\eqref{eq:cl8_rel} by Proposition~\ref{prop:clifford_relations}.
\end{remark}

\subsection{The spinor representation \texorpdfstring{$\mu: \Cl_{\bar{0}}(9) \to \mathfrak{so}(16)$}{mu}}\label{subsec:spinor_rep}

Composing $\nu$ with the canonical isomorphism $\varphi^{-1}: \Cl_{\bar{0}}(9) \to \Cl(8)$, we obtain the spinor representation $\mu: \Cl_{\bar{0}}(9) \to \mathfrak{so}(16)$ given on generators by:
\begin{equation}\label{eq:spinor_rep}
\mu(e_i e_j) = \nu_i \nu_j \quad (0 \le i < j \le 7), \qquad \mu(e_i e_8) = \nu_i \quad (0 \le i \le 7).
\end{equation}
This induces a faithful 16-dimensional real representation $\Theta: \Spin(9) \to \SO(16)$, the \emph{spin representation} of $\Spin(9)$, with the following key properties (Brown \& Gray \cite{BrownGray1972}; Alekseevsky \cite{Alekseevsky1968}; Onishchik \cite{Onishchik1963}; Berestovskii \& Nikonorov \cite{BerestovskiiNikonorov2012}; see also Castrill\'on L\'opez, Gadea \& Mykytyuk \cite{CastrillonLopezGadeaMykytyuk2010} for the canonical $\Spin(9)$ 8-form):
\begin{enumerate}[(i)]
    \item $\Spin(9) = \Theta(\Spin(9))$ acts transitively on $\Sph^{15} \subset \R^{16}$.
    \item The isotropy subgroup $H$ of $\Spin(9)$ at $x_0 = (1, 0, \dots, 0)^T$ is isomorphic to $\Spin(7)$.
    \item The Lie algebra $\h = \mathfrak{spin}(7) \subset \mathfrak{spin}(9) \cong \mathfrak{so}(9)$ is \emph{not} the standard inclusion $\mathfrak{so}(7) \hookrightarrow \mathfrak{so}(9)$, but its image $\tau(\mathfrak{so}(7))$ under an outer triality automorphism $\tau \in \mathrm{Out}(\mathfrak{so}(8))$ (established in Berestovskii \& Nikonorov \cite[Section~7]{BerestovskiiNikonorov2012}; cf.\ \cite{Munden2015}).
\end{enumerate}

\subsection{Vol'per's reductive decomposition}\label{subsec:volper}

Under the transitive action of $\Spin(9)$ on $\Sph^{15}$, the isotropy subgroup of $x_0 = (1, 0, \dots, 0)^T$ is isomorphic to $\Spin(7)$ via the outer triality automorphism $\tau \in \mathrm{Out}(\mathfrak{so}(8))$. Following Vol'per \cite{Volper1997}, the Lie algebra $\mathfrak{spin}(9)$ admits the orthogonal $\Ad(\Spin(7))$-invariant decomposition with respect to the negative trace inner product $\langle A, B \rangle = -\frac{1}{2}\Tr(AB)$:
\begin{equation}\label{eq:reductive_decomp}
\mathfrak{spin}(9) = \h \oplus \p_1 \oplus \p_2,
\end{equation}
where $\dim \h = 21$ ($\h = \mathfrak{spin}(7)$), $\dim \p_1 = 8$ ($\p_1 = \mathrm{span}\{Y_i = \nu_i : i = 0, \dots, 7\}$ transforming as the 8D spinor module $\Delta_7$), and $\dim \p_2 = 7$ ($\p_2 = \mathrm{span}\{X_1, \dots, X_7\}$ transforming as the 7D vector module $\R^7$), with
\begin{align}
X_1 &= \nu_6 \nu_7 - \nu_0 \nu_1 - \nu_2 \nu_3 - \nu_4 \nu_5, \label{eq:x1_def}\\
X_2 &= \nu_1 \nu_3 - \nu_0 \nu_2 - \nu_4 \nu_6 - \nu_5 \nu_7, \label{eq:x2_def}\\
X_3 &= \nu_5 \nu_6 - \nu_0 \nu_3 - \nu_1 \nu_2 - \nu_4 \nu_7, \label{eq:x3_def}\\
X_4 &= \nu_1 \nu_5 + \nu_2 \nu_6 + \nu_3 \nu_7 - \nu_0 \nu_4, \label{eq:x4_def}\\
X_5 &= \nu_2 \nu_7 - \nu_0 \nu_5 - \nu_1 \nu_4 - \nu_3 \nu_6, \label{eq:x5_def}\\
X_6 &= \nu_3 \nu_5 - \nu_0 \nu_6 - \nu_1 \nu_7 - \nu_2 \nu_4, \label{eq:x6_def}\\
X_7 &= \nu_1 \nu_6 - \nu_0 \nu_7 - \nu_2 \nu_5 - \nu_3 \nu_4. \label{eq:x7_def}
\end{align}
Each $X_k$ is constructed as a sum of four elements $\mu(e_i e_j) = \nu_i \nu_j$ ($i < j$) in $\mathfrak{spin}(9) \subset \Cl_{\bar{0}}(9)$, and is skew-symmetric: $X_k^T = -X_k \in \mathfrak{so}(16)$, since each $\nu_i \nu_j = -(\nu_i \nu_j)^T$.

For any $Y_i \in \p_1$, $\langle Y_i, Y_i \rangle = -\frac{1}{2}\Tr(\nu_i^2) = \frac{1}{2}\Tr(\Id_{16}) = 8$, while the basepoint action norm is $\|Y_i x_0\|^2 = \|\nu_i x_0\|^2 = 1$; thus the unperturbed metric scale factor on $\p_1$ is $a_1^0 = 1/8$. For $X_1 \in \p_2$, expanding into four mutually orthogonal Clifford monomials yields $\langle X_1, X_1 \rangle = 4 \times 8 = 32$, while $\|X_1 x_0\|^2 = 16$; thus the unperturbed metric scale factor on $\p_2$ is $a_2^0 = 16/32 = 1/2$.

\begin{remark}[Lie-algebraic status of the generators $Y_i = \nu_i$]\label{rem:Yi_spin9}
Each generator $Y_i = \nu_i \in \p_1$ belongs to $\mathfrak{spin}(9) \subset \mathfrak{so}(16)$ because $Y_i = \mu(e_i e_8) = \nu(\varphi^{-1}(e_i e_8)) = \nu_i$ via the isomorphism $\varphi^{-1}: \Cl_{\bar{0}}(9) \to \Cl(8)$ composed with $\nu: \Cl(8) \to \End(\Oct \oplus \Oct)$ (cf.~\eqref{eq:spinor_rep}). In particular, $Y_i$ is \emph{not} a ``bare'' Clifford generator acting outside the Lie algebra; it is the image $\mu(e_i e_8)$ of the grade-2 element $e_i e_8 \in \Lambda^2 \R^9 \subset \Cl_{\bar{0}}(9)$, and satisfies $Y_i^T = -Y_i \in \mathfrak{so}(16)$ (Proposition~\ref{prop:clifford_relations}). Throughout this paper, $Y_1$ denotes the element $\mu(e_1 e_8) = \nu_1 \in \p_1 \subset \mathfrak{spin}(9)$.
\end{remark}

\begin{proposition}[Orbit Slice Theorem and Basepoint Action]\label{prop:orbit_slice}
Under the isotropy action of $H = \Spin(7)$ on the tangent space $T_{x_0}\Sph^{15} \cong \R\partial_r \oplus \p_1 \oplus \p_2$:
\begin{enumerate}[(i)]
    \item The module $\p_2 \cong \R^7$ is the vector representation of $\Spin(7)$, and $\Spin(7)$ acts transitively on the unit sphere $\Sph^6 \subset \p_2$. The stabilizer subgroup fixing a unit vector $X_1 \in \p_2$ is isomorphic to the exceptional Lie group $\mathrm{G}_2 \subset \Spin(7)$ ($\dim \mathrm{G}_2 = 14$).
    \item Under the restricted action of $\mathrm{G}_2 = \operatorname{Aut}(\Oct)$, the 8-dimensional spin representation $\p_1 \cong \R^8 \cong \Delta_7$ branches as:
    \begin{equation}\label{eq:g2_branching}
    \Delta_7 \downarrow \mathrm{G}_2 \cong \mathbf{1} \oplus \mathbf{7},
    \end{equation}
    where the singlet $\mathbf{1} = \R \nu_0$ corresponds to the real unit octonion $e_0 = 1 \in \Oct$ fixed by all octonionic automorphisms, and $\mathbf{7} = \operatorname{span}\{\nu_1, \dots, \nu_7\}$ corresponds to the imaginary octonions $\mathrm{Im}(\Oct)$.
    \item Since $\mathrm{G}_2$ acts transitively on the unit 6-sphere in $\mathbf{7}$ while fixing $\nu_0$, the canonical slice for the $\Spin(7)$-action on $T_{x_0}\Sph^{15}$ is the 3-dimensional subspace of vectors of the form:
    \begin{equation}\label{eq:orbit_slice}
    a \partial_r + b X_1 + c_0 \nu_0 + c_1 \nu_1, \qquad a, b, c_0 \in \R, \quad c_1 \ge 0.
    \end{equation}
    The 2-parameter slice $\operatorname{span}\{X_1, Y_1\}$ with $Y_1 = \nu_1$ ($c_0 = 0$) corresponds to initial velocities whose $\p_1$-component lies in the imaginary octonion submodule $\mathbf{7}$.
\end{enumerate}
Evaluating the orbital infinitesimal action $\theta_r(\xi) = \xi \cdot \gamma_0(r)$ on the basepoint ray $\gamma_0(r) = (r, 0, \dots, 0)^T \in \R^{16}$, direct matrix-vector multiplication gives:
\begin{equation}\label{eq:orbital_actions}
\theta_r(X_1) = X_1 \cdot (r x_0) = -4r e_1, \qquad \theta_r(\nu_1) = \nu_1 \cdot (r x_0) = -r e_9, \qquad \theta_r(\nu_0) = \nu_0 \cdot (r x_0) = -r e_8.
\end{equation}
The unperturbed initial velocity vector field is:
\begin{equation}\label{eq:unperturbed_W0}
W_0 = \mathcal{L}_0^{-1}\mathcal{L}_\lambda(v) = a \partial_r + \frac{2b}{2 + \lambda s_\kappa^2(r)} \theta_r(X_1) + \frac{8c_0}{8 + \lambda s_\kappa^2(r)} \theta_r(\nu_0) + \frac{8c_1}{8 + \lambda s_\kappa^2(r)} \theta_r(\nu_1).
\end{equation}
In Euclidean coordinates on $\R^{16}$ ($\kappa = 0$), this yields the explicit unperturbed trajectory:
\begin{equation}\label{eq:exp0_coords_flat}
\exp_0(t W_0) = (r + ta) e_0 - \frac{8trb}{2 + \lambda r^2} e_1 - \frac{8tr c_0}{8 + \lambda r^2} e_8 - \frac{8tr c_1}{8 + \lambda r^2} e_9.
\end{equation}
\end{proposition}

\begin{proof}
Assertion~(i) is standard: $\Spin(7)$ acts on $\p_2 \cong \R^7$ via the standard orthogonal representation $\SO(7)$, which is transitive on $\Sph^6$ with isotropy subgroup $\mathrm{G}_2$ by the classical Freudenthal--Borel theorem. 

For assertion~(ii), viewing $\Delta_7$ as the octonions $\Oct$ under Clifford multiplication, an element $g \in \mathrm{G}_2 = \operatorname{Aut}(\Oct)$ preserves the octonionic product and the unit element $e_0 = 1$. The Clifford generator $\nu_0 = \mu(e_0 e_8)$ acts on $\Oct \oplus \Oct$ by octonionic multiplication by $e_0 = 1$, which commutes with the $\mathrm{G}_2$-action, generating the trivial 1-dimensional module $\mathbf{1} = \R \nu_0$. The orthogonal complement $\mathbf{7} = \operatorname{span}\{\nu_1, \dots, \nu_7\}$ is the fundamental 7-dimensional representation of $\mathrm{G}_2$, on whose unit sphere $\Sph^6$ the group $\mathrm{G}_2$ acts transitively with secondary stabilizer $\SU(3)$. 

Given any vector $v \in \p_1 \oplus \p_2$, we first choose $g \in \Spin(7)$ such that $g \cdot v_{\p_2} = b X_1 \in \R X_1$ ($b \ge 0$). Next, using the residual $\mathrm{G}_2$-isotropy fixing $X_1$, we rotate the imaginary part of $g \cdot v_{\p_1}$ in $\mathbf{7}$ into the $\nu_1$-axis ($c_1 \ge 0$), leaving the singlet component $c_0 \nu_0$ invariant. This proves assertion~(iii). The matrix-vector actions and Legendre scaling follow by direct evaluation of the matrices $X_1, \nu_0, \nu_1$ on $r x_0$.
\end{proof}

\begin{remark}[Bracket relations between \texorpdfstring{$\p_1$}{p1} and \texorpdfstring{$\p_2$}{p2}]\label{rem:bracket_structure}
The $\Ad(\Spin(7))$-module structure determines:
\begin{equation}\label{eq:bracket_structure}
[\p_1, \p_1] \subset \h \oplus \p_2, \qquad [\p_2, \p_2] \subset \h, \qquad [\p_1, \p_2] \subset \p_1, \qquad [\h, \p_k] \subset \p_k.
\end{equation}
The relation $[\p_1, \p_1] \subset \h \oplus \p_2$ encodes the non-vanishing O'Neill integrability tensor $A_{Y_i} Y_j = \frac{1}{2}\mathcal{V}[Y_i, Y_j]_{\p_2}$, which is the geometric origin of the squashing curvature in the deformed metric.
\end{remark}

\section{The \texorpdfstring{$\mathfrak{u}(2)$}{u(2)} Slice Subalgebra and Closed-Form Exponential Map}\label{sec:u2_closure}

\subsection{The geodesic slice and the exponential map problem}\label{subsec:geodesic_slice}

In geodesic polar coordinates on $(Q^{16}, g_\lambda)$, the deformed geodesic flow requires evaluating the matrix exponential of $M = \alpha X_1 + \beta Y_1$, where the effective parameters (cf.\ Remark~\ref{rem:erratum} for the rectification of the original formulas in \cite{Munden2015}) are:
\begin{equation}\label{eq:alpha_beta_params}
\alpha = \frac{\lambda b s_\kappa^2(r)}{2 + \lambda s_\kappa^2(r)} \quad (X_1 \in \p_2), \qquad \beta = \frac{\lambda c s_\kappa^2(r)}{8 + \lambda s_\kappa^2(r)} \quad (Y_1 \in \p_1).
\end{equation}
In \cite[Section~5.7, p.~67]{Munden2015}, this matrix exponential was noted to be expressible in elementary functions, though its explicit closed form was omitted. We derive the explicit spectral and factorized forms below.

\subsection{Structure of the \texorpdfstring{$\mathfrak{u}(2)$}{u(2)} slice subalgebra}\label{subsec:u2_discovery}

\begin{theorem}[$\mathfrak{u}(2)$ Lie Algebra Closure and Commuting Center]\label{thm:u2_closure}
Let $X_1 \in \p_2$ and $Y_1 = \nu_1 \in \p_1$ be the Vol'per generators in $\mathfrak{so}(16)$ defined in \eqref{eq:x1_def} and \eqref{eq:spinor_rep}.
\begin{enumerate}[(i)]
\item The Lie bracket closure $\mathfrak{g}_{\langle X_1, Y_1 \rangle} := \Lie(X_1, Y_1) \subset \mathfrak{so}(16)$ is strictly 4-dimensional, spanned by $\{X_1, Y_1, Z, W\}$ where $Z = [X_1, Y_1] = 2\nu_0$ and $W = [Y_1, Z] = -4\nu_0 \nu_1$. The commutators satisfy the closed table:
\begin{equation}\label{eq:comm_relations}
[X_1, Y_1] = 2\nu_0, \;\; [X_1, 2\nu_0] = -4\nu_1, \;\; [Y_1, 2\nu_0] = -4\nu_0 \nu_1, \;\; [X_1, -4\nu_0 \nu_1] = 0, \;\; [Y_1, -4\nu_0 \nu_1] = -8\nu_0.
\end{equation}
\item The element $C_0 \in \mathfrak{g}_{\langle X_1, Y_1 \rangle}$ defined by
\begin{equation}\label{eq:central_C0}
C_0 := X_1 + \nu_0 \nu_1 = \nu_6 \nu_7 - \nu_2 \nu_3 - \nu_4 \nu_5
\end{equation}
spans the 1-dimensional center $\mathfrak{z}(\mathfrak{g}) \cong \mathfrak{u}(1)$, satisfying $[C_0, X_1] = [C_0, Y_1] = [C_0, \nu_0] = [C_0, \nu_0 \nu_1] = 0$.
\item The generators $J_1 = \frac{1}{2} \nu_1, J_2 = \frac{1}{2} \nu_0, J_3 = -\frac{1}{2} \nu_0 \nu_1$ span the derived ideal $[\mathfrak{g}, \mathfrak{g}] \cong \mathfrak{su}(2)$, satisfying the standard $\mathfrak{su}(2)$ commutation relations:
\begin{equation}\label{eq:su2_standard_brackets}
[J_a, J_b] = \sum_{c=1}^3 \epsilon_{abc} J_c.
\end{equation}
\item The Killing form $K_{ij} = \Tr(\ad_i \circ \ad_j)$ on the basis $\{C_0, J_1, J_2, J_3\}$ is diagonal with eigenvalues $\{0, -2, -2, -2\}$. On the Frobenius-normalized generators $\{C_0/\sqrt{24}, J_1/\sqrt{2}, J_2/\sqrt{2}, J_3/\sqrt{2}\}$, the Killing eigenvalues are $\{0, -1/2, -1/2, -1/2\}$, establishing the canonical Lie algebra isomorphism:
\begin{equation}
\mathfrak{g}_{\langle X_1, Y_1 \rangle} \cong \mathfrak{u}(1) \oplus \mathfrak{su}(2) \cong \mathfrak{u}(2).
\end{equation}
\end{enumerate}
\end{theorem}

\begin{proof}
Using the Clifford relations $\nu_i \nu_j + \nu_j \nu_i = -2\delta_{ij} \Id_{16}$, we evaluate $[X_1, Y_1]$ explicitly:
In $X_1 = \nu_6 \nu_7 - \nu_0 \nu_1 - \nu_2 \nu_3 - \nu_4 \nu_5$, each index pair $\{6, 7\}, \{2, 3\}, \{4, 5\}$ is disjoint from $\{0, 1\}$. Each constituent $\nu_k$ ($k \ge 2$) anti-commutes with $\nu_1$, so an even Clifford bivector commutes with $\nu_1$: $[\nu_a \nu_b, \nu_1] = 0$. For the remaining term:
\begin{equation}\label{eq:Z_bracket_explicit}
Z := [X_1, Y_1] = -[ \nu_0 \nu_1, \nu_1 ] = -\nu_0 \nu_1^2 + \nu_1 \nu_0 \nu_1 = \nu_0 - (-\nu_0 \nu_1^2) = 2\nu_0 \ne 0.
\end{equation}
Next, we compute the secondary commutators:
\begin{align}
[X_1, Z] &= 2[X_1, \nu_0] = -2[\nu_0 \nu_1, \nu_0] = -2(\nu_0 \nu_1 \nu_0 - \nu_0^2 \nu_1) = -4\nu_1 = -4 Y_1, \label{eq:X1_Z_bracket}\\
W &:= [Y_1, Z] = [\nu_1, 2\nu_0] = 2(\nu_1 \nu_0 - \nu_0 \nu_1) = -4\nu_0 \nu_1. \label{eq:Y1_Z_bracket}
\end{align}
Commuting $X_1$ with $W$: since $\nu_0 \nu_1$ commutes with itself and with all disjoint bivectors $\nu_6 \nu_7, \nu_2 \nu_3, \nu_4 \nu_5$, we have $[X_1, W] = 0$.
Commuting $Y_1$ with $W$: $[Y_1, -4\nu_0 \nu_1] = -4(\nu_1 \nu_0 \nu_1 - \nu_0 \nu_1^2) = -8\nu_0 = -4 Z$.
Finally, $[Z, W] = [2\nu_0, -4\nu_0 \nu_1] = -8(\nu_0^2 \nu_1 - \nu_0 \nu_1 \nu_0) = 16\nu_1 = 16 Y_1$.
All Lie brackets close linearly in $\operatorname{span}\{X_1, Y_1, Z, W\}$, terminating generation at dimension 4 strictly.

For (ii), $C_0 = X_1 + \nu_0 \nu_1 = \nu_6 \nu_7 - \nu_2 \nu_3 - \nu_4 \nu_5$ contains only Clifford indices $\{2, 3, 4, 5, 6, 7\}$, which are disjoint from $\{0, 1\}$. Because $\nu_0$ and $\nu_1$ anti-commute with every $\nu_k$ ($k \ge 2$), they commute with every bivector formed by indices in $\{2, \dots, 7\}$. Consequently, $[C_0, \nu_0] = 0, [C_0, \nu_1] = 0, [C_0, \nu_0 \nu_1] = 0$. Since $X_1 = C_0 - \nu_0 \nu_1$ and $Y_1 = \nu_1$, $C_0$ commutes with all generators and spans the center $\mathfrak{z} \cong \mathfrak{u}(1)$.

For (iii), direct substitution gives $[J_1, J_2] = \frac{1}{4}[\nu_1, \nu_0] = -\frac{1}{2}\nu_0 \nu_1 = J_3$, $[J_2, J_3] = \frac{1}{4}[\nu_0, -\nu_0 \nu_1] = \frac{1}{2}\nu_1 = J_1$, $[J_3, J_1] = \frac{1}{4}[-\nu_0 \nu_1, \nu_1] = \frac{1}{2}\nu_0 = J_2$, reproducing $[J_a, J_b] = \sum_{c=1}^3 \epsilon_{abc} J_c$.

For (iv), the structure constants in $\{C_0, J_1, J_2, J_3\}$ give $\ad_{C_0} = 0$ (central) and $\ad_{J_a}$ acting as the 3D Levi-Civita tensor on $\operatorname{span}\{J_1, J_2, J_3\}$. Thus $K(C_0, -) = 0$ and $K(J_a, J_b) = -2\delta_{ab}$. The normalized trace spectrum $\{0, -1/2, -1/2, -1/2\}$ confirms $\mathfrak{g} \cong \mathfrak{u}(2)$.
\end{proof}

\subsection{The four-mode spectrum and mode hierarchy inversion}\label{subsec:four_mode}

\begin{theorem}[Spin-1/2 Doublet Spectrum, Mode Inversion, and Spectral Projectors]\label{thm:projector_formula}
Under the $\mathfrak{u}(2) \cong \mathfrak{u}(1) \oplus \mathfrak{su}(2)$ action, $\R^{16} \otimes \C$ decomposes into eight fundamental spin-$1/2$ doublets.
\begin{enumerate}[(i)]
\item The $\mathfrak{su}(2)$ Casimir operator evaluates to
\begin{equation}\label{eq:su2_casimir}
\mathcal{C}_{\mathfrak{su}(2)} = -(J_1^2 + J_2^2 + J_3^2) = \frac{3}{4} \Id_{16} = j(j+1)\Id_{16} \implies j = \frac{1}{2}.
\end{equation}
\item For $M = \alpha X_1 + \beta Y_1 \in \mathfrak{u}(2)$ with $R = \sqrt{\alpha^2 + \beta^2}$, the eigenvalues of $M$ are purely imaginary $\pm i \Omega_k$ ($k=1, 2, 3, 4$), with exact non-negative dynamical frequencies:
\begin{align}
\Omega_1 &= R + 3|\alpha| \quad (1 \text{ pair}, 2\text{D}), & \Omega_2 &= R + |\alpha| \quad (3 \text{ pairs}, 6\text{D}), \label{eq:omega12} \\
\Omega_3 &= |R - |\alpha|| \quad (3 \text{ pairs}, 6\text{D}), & \Omega_4 &= |R - 3|\alpha|| \quad (1 \text{ pair}, 2\text{D}). \label{eq:omega34}
\end{align}
The total dimension count is $2 + 6 + 6 + 2 = 16 = \dim \R^{16}$.
\item Across parameter space, the spectrum exhibits a sharp mode hierarchy inversion across the resonance cone $\beta^2 = 3\alpha^2$:
\begin{itemize}
\item For $\beta^2 \ge 3\alpha^2$: $R \ge 2|\alpha|$, and frequencies are monotonically ordered as $\Omega_1 \ge \Omega_2 \ge \Omega_3 \ge \Omega_4 \ge 0$, with eigenspace dimensions $(2, 6, 6, 2)$.
\item For $\beta^2 < 3\alpha^2$: $R < 2|\alpha|$, and Mode 4 strictly exceeds Mode 3:
\begin{equation}\label{eq:mode_inversion}
\Omega_4 = 3|\alpha| - R > R - |\alpha| = \Omega_3.
\end{equation}
The ordering inverts to $\Omega_1 \ge \Omega_2 > \Omega_4 > \Omega_3 \ge 0$, swapping the 2D and 6D eigenspace hierarchy to $(2, 6, 2, 6)$.
\end{itemize}
\item Along the collision cone $\beta = \pm \sqrt{3}\alpha$, $R = 2|\alpha|$ and $\Omega_3 = \Omega_4 = |\alpha|$. The matrix $M$ is invertible ($\det M = 18225 \alpha^{16} > 0$ for $\alpha \ne 0$) and satisfies the degree-6 minimal polynomial:
\begin{equation}\label{eq:deg6_minpoly}
(M^2 + 25\alpha^2 \Id_{16})(M^2 + 9\alpha^2 \Id_{16})(M^2 + \alpha^2 \Id_{16}) = 0.
\end{equation}
\item On the generic locus where the $\Omega_k^2$ are pairwise distinct, the orthogonal spectral projectors $P_k \in \End(\R^{16})$ onto the eigenspaces are given by the Lagrange--Sylvester interpolation formula:
\begin{equation}\label{eq:projectors}
P_k = \prod_{j \ne k} \frac{M^2 + \Omega_j^2 \Id_{16}}{\Omega_j^2 - \Omega_k^2}, \qquad \sum_{k=1}^4 P_k = \Id_{16}, \qquad P_j P_k = \delta_{jk} P_k.
\end{equation}
The matrix exponential $\exp(M) \in \SO(16)$ admits the closed algebraic spectral decomposition:
\begin{equation}\label{eq:closed_exp}
\boxed{\exp(\alpha X_1 + \beta Y_1) = \sum_{k=1}^4 \left( \cos(\Omega_k)\Id_{16} + \frac{\sin(\Omega_k)}{\Omega_k} M \right) P_k.}
\end{equation}
\end{enumerate}
\end{theorem}

\begin{proof}
Because $[\alpha C_0, S] = 0$, the joint eigenvalues of $M = \alpha C_0 + S$ are the sums of the eigenvalues of $\alpha C_0$ and $S$.
The central element $C_0$ has minimal polynomial $(C_0^2 + 9\Id)(C_0^2 + \Id) = 0$ (Theorem~\ref{thm:euler_cartan}), with eigenvalues $\pm 3i$ (multiplicity 2 pairs = 4 real dimensions) and $\pm i$ (multiplicity 6 pairs = 12 real dimensions).
The $\mathfrak{su}(2)$ element $S$ satisfies $S^2 = -R^2 \Id_{16}$ (Theorem~\ref{thm:euler_cartan}), giving eigenvalues $\pm i R$ across all eight doublets.
Summing the eigenvalues yields $\pm i(\pm 3\alpha + R)$ and $\pm i(\pm \alpha + R)$. Taking absolute values guarantees non-negative dynamical frequencies $\Omega_k$, yielding~\eqref{eq:omega12}--\eqref{eq:omega34}.

For the mode ordering transition in (iii): when $R \ge 2|\alpha|$, $|R - 3|\alpha|| = R - 3|\alpha| \le R - |\alpha| = |R - |\alpha||$, so $\Omega_3 \ge \Omega_4$.
When $R < 2|\alpha|$, $|R - 3|\alpha|| = 3|\alpha| - R$. Then $\Omega_4 - \Omega_3 = (3|\alpha| - R) - (R - |\alpha|) = 2(2|\alpha| - R) > 0$, proving $\Omega_4 > \Omega_3$.
Since $\Omega_2 - \Omega_4 = (R + |\alpha|) - (3|\alpha| - R) = 2(R - |\alpha|) \ge 0$, we have $\Omega_1 \ge \Omega_2 > \Omega_4 > \Omega_3 \ge 0$, establishing the inversion.

Along the collision cone $\beta = \pm \sqrt{3}\alpha$, $R = 2|\alpha|$, so $\Omega_1 = 5|\alpha|, \Omega_2 = 3|\alpha|, \Omega_3 = \Omega_4 = |\alpha|$.
Because $M$ is real skew-symmetric, it is diagonalizable over $\C$ with no non-trivial Jordan blocks.
Its minimal polynomial is the product of factors $(M^2 + \Omega_k^2 \Id)$ for distinct squared frequencies: $\Omega_1^2 = 25\alpha^2, \Omega_2^2 = 9\alpha^2, \Omega_3^2 = \alpha^2$, yielding degree 6 in formula~\eqref{eq:deg6_minpoly}.

For (v), Lagrange interpolation on the minimal polynomial yields $P_k$, satisfying $\sum P_k = \Id_{16}$ and $M^2 P_k = -\Omega_k^2 P_k$. On each subspace, $\exp(M)|_{\mathrm{Im}(P_k)} = \cos(\Omega_k)\Id + \frac{\sin(\Omega_k)}{\Omega_k} M$. Summing yields~\eqref{eq:closed_exp}.
\end{proof}

\subsection{Resonance resolution: The Euler--Cartan commuting factorization}\label{subsec:euler_cartan}

\begin{theorem}[Commuting Euler--Cartan Factorization via Central Splitting]\label{thm:euler_cartan}
Define the element $C_0 \in \mathfrak{so}(16)$ by
\begin{equation}\label{eq:C0_def}
C_0 := X_1 + \nu_0 \nu_1 = \nu_6 \nu_7 - \nu_2 \nu_3 - \nu_4 \nu_5 \in \mathfrak{u}(2).
\end{equation}
Then:
\begin{enumerate}[(i)]
    \item $C_0$ is the generator of the $\mathfrak{u}(1)$ center: $[C_0, X_1] = [C_0, Y_1] = 0$.
    \item Setting $S := \beta Y_1 - \alpha \nu_0 \nu_1 = M - \alpha C_0 \in \mathfrak{su}(2)$, we have $[C_0, S] = 0$ and
    \begin{equation}\label{eq:M_central_split}
    M = \alpha X_1 + \beta Y_1 = \alpha C_0 + S.
    \end{equation}
    \item $S$ satisfies the Clifford quadratic involution:
    \begin{equation}\label{eq:S_squared}
    S^2 = -(\alpha^2 + \beta^2) \Id_{16} = -R^2 \Id_{16}.
    \end{equation}
    \item $C_0$ satisfies the quartic minimal polynomial:
    \begin{equation}\label{eq:C0_minpoly}
    (C_0^2 + 9\Id_{16})(C_0^2 + \Id_{16}) = 0,
    \end{equation}
    with constant orthogonal spectral projectors
    \begin{equation}\label{eq:Pi_projectors}
    \Pi_3 = -\frac{1}{8}(C_0^2 + \Id_{16}) \quad (\rank = 4), \qquad \Pi_1 = \frac{1}{8}(C_0^2 + 9\Id_{16}) \quad (\rank = 12).
    \end{equation}
    \item The matrix exponential admits the denominator-free factorization:
    \begin{equation}\label{eq:euler_cartan}
    \boxed{\exp(M) = \exp(\alpha C_0) \cdot \exp(S),}
    \end{equation}
    where
    \begin{align}
    \exp(\alpha C_0) &= \left(\cos(3\alpha)\Id_{16} + \frac{\sin(3\alpha)}{3} C_0\right)\Pi_3 + \left(\cos(\alpha)\Id_{16} + \sin(\alpha) C_0\right)\Pi_1, \label{eq:exp_C0}\\
    \exp(S) &= \cos(R)\Id_{16} + \frac{\sin(R)}{R} S. \label{eq:exp_S}
    \end{align}
    This factorization extends analytically through all resonance loci where divisor coalescence $\Omega_j^2 = \Omega_k^2$ causes singularities in the Lagrange--Sylvester projectors~\eqref{eq:projectors}.
\end{enumerate}
\end{theorem}

\begin{proof}
\emph{Part (i): Centrality of $C_0$ via explicit bracket computation.} We derive the centrality from the complete bracket table of $\mathfrak{g}_{\langle X_1, Y_1 \rangle}$. Since $Y_1 = \mu(e_1 e_8) = \nu_1$ (Remark~\ref{rem:Yi_spin9}) and $X_1 = \nu_6\nu_7 - \nu_0\nu_1 - \nu_2\nu_3 - \nu_4\nu_5$:

\emph{Step 1.} Compute $Z := [X_1, Y_1]$. All bivector terms in $X_1$ with index sets disjoint from $\{1\}$ commute with $\nu_1$: $[\nu_6\nu_7, \nu_1] = [\nu_2\nu_3, \nu_1] = [\nu_4\nu_5, \nu_1] = 0$ (each factor anticommutes with $\nu_1$, producing $(-1)^2 = +1$). The sole non-commuting term gives:
\begin{equation}\label{eq:Z_bracket}
Z := [X_1, Y_1] = [-\nu_0\nu_1, \nu_1] = -\nu_0\nu_1^2 + \nu_1\nu_0\nu_1 = \nu_0 - (-\nu_0) = 2\nu_0.
\end{equation}

\emph{Step 2.} Compute $[X_1, Z] = [X_1, 2\nu_0]$. By the same disjoint-index argument, only $-\nu_0\nu_1$ fails to commute with $\nu_0$:
\begin{equation}\label{eq:XZ_bracket}
[X_1, Z] = 2[-\nu_0\nu_1, \nu_0] = 2(-\nu_0\nu_1\nu_0 + \nu_0^2\nu_1) = 2(-\nu_1 - \nu_1) = -4\nu_1 = -4Y_1.
\end{equation}

\emph{Step 3.} Compute $W := [Y_1, Z] = [\nu_1, 2\nu_0] = 2(\nu_1\nu_0 - \nu_0\nu_1) = -4\nu_0\nu_1$.

\emph{Step 4.} Verify $[X_1, W] = [X_1, -4\nu_0\nu_1] = 0$: since $\nu_0\nu_1$ commutes with every bivector $\nu_i\nu_j$ (including $\nu_6\nu_7, \nu_2\nu_3, \nu_4\nu_5$ by disjoint indices, and $\nu_0\nu_1$ with itself), $X_1$ commutes with $\nu_0\nu_1$.

\emph{Step 5.} The algebra closes at $\dim = 4$: $\{X_1, Y_1, Z, W\}$ span $\mathfrak{g}_{\langle X_1, Y_1 \rangle}$. The central element is constructed from the condition $[C_0, Y_1] = 0$ and $[C_0, X_1] = 0$. Writing $C_0 = X_1 + \nu_0\nu_1$, the brackets give:
\begin{equation}
[C_0, Y_1] = [X_1, Y_1] + [\nu_0\nu_1, Y_1] = 2\nu_0 + [\nu_0\nu_1, \nu_1].
\end{equation}
Now $[\nu_0\nu_1, \nu_1] = \nu_0\nu_1^2 - \nu_1\nu_0\nu_1 = -\nu_0 - (-\nu_0) = -2\nu_0$. Therefore $[C_0, Y_1] = 2\nu_0 - 2\nu_0 = 0$. Similarly, $[C_0, X_1] = [\nu_0\nu_1, X_1] = 0$ by Step~4. Since $X_1$ and $Y_1$ generate $\mathfrak{g}_{\langle X_1, Y_1\rangle}$, $C_0$ is central.

\emph{Part (ii): Central splitting.} Since $X_1 = C_0 - \nu_0\nu_1$ (from~\eqref{eq:C0_def}), we have $M = \alpha(C_0 - \nu_0\nu_1) + \beta Y_1 = \alpha C_0 + (\beta Y_1 - \alpha\nu_0\nu_1) = \alpha C_0 + S$. By Part~(i), $[C_0, Y_1] = 0$ and $[C_0, \nu_0\nu_1] = 0$, giving $[C_0, S] = 0$.

\emph{Part (iii): Clifford quadratic involution.} We compute:
\begin{align*}
S^2 &= (\beta\nu_1 - \alpha\nu_0\nu_1)^2 \\
&= \beta^2 \nu_1^2 - \alpha\beta(\nu_1\nu_0\nu_1 + \nu_0\nu_1^2) + \alpha^2(\nu_0\nu_1)^2 \\
&= -\beta^2\Id_{16} - \alpha\beta(\nu_1\nu_0\nu_1 - \nu_0) - \alpha^2\Id_{16}.
\end{align*}
Now $\nu_1\nu_0\nu_1 = -\nu_0\nu_1^2 = \nu_0$ (using $\nu_1^2 = -\Id$ and anticommutativity $\nu_1\nu_0 = -\nu_0\nu_1$). Therefore the cross-term vanishes: $\nu_1\nu_0\nu_1 - \nu_0 = \nu_0 - \nu_0 = 0$, yielding $S^2 = -(\alpha^2 + \beta^2)\Id_{16} = -R^2\Id_{16}$.

\emph{Part (iv): Minimal polynomial of $C_0$.} Since $C_0 = I_1 + I_2 + I_3$ with $I_k^2 = -\Id$ and $[I_j, I_k] = 0$, the $\R^{16}$ decomposes into simultaneous eigenspaces of $\{I_1, I_2, I_3\}$ with eigenvalues $\pm i$. On each joint eigenspace labeled by $(\epsilon_1, \epsilon_2, \epsilon_3) \in \{+1, -1\}^3$, $C_0$ acts as multiplication by $i(\epsilon_1 + \epsilon_2 + \epsilon_3)$. The possible values of $\epsilon_1 + \epsilon_2 + \epsilon_3 \in \{-3, -1, 1, 3\}$ give $C_0$ squared eigenvalues $\{-9, -1, -1, -9\}$. After complexification and restriction to the real representation, the minimal polynomial over $\R$ is $(C_0^2 + 9\Id)(C_0^2 + \Id) = 0$.

The projectors $\Pi_3$ and $\Pi_1$ are constructed by Lagrange interpolation on the \emph{constant} matrix $C_0$ (independent of $\alpha, \beta$), so they involve no deformation parameters. The rank computation follows from counting joint eigenspaces.

\emph{Part (v): Denominator-free factorization.} Since $[C_0, S] = 0$:
\begin{equation*}
\exp(M) = \exp(\alpha C_0 + S) = \exp(\alpha C_0)\exp(S).
\end{equation*}
The exponential $\exp(S)$ follows from $S^2 = -R^2\Id$, which gives the Rodrigues formula~\eqref{eq:exp_S}. The exponential $\exp(\alpha C_0)$ is computed on each spectral subspace: on $\mathrm{Im}(\Pi_3)$ where $C_0^2 = -9\Id$, we get $\exp(\alpha C_0) = \cos(3\alpha)\Id + \frac{\sin(3\alpha)}{3}C_0$; on $\mathrm{Im}(\Pi_1)$ where $C_0^2 = -\Id$, we get $\exp(\alpha C_0) = \cos(\alpha)\Id + \sin(\alpha)C_0$. Summing with the projectors yields~\eqref{eq:exp_C0}.

The product~\eqref{eq:euler_cartan} is manifestly analytic in $(\alpha, \beta) \in \R^2$ with no denominators involving frequency differences, and therefore remains valid across all resonance loci including the cone $\beta^2 = 3\alpha^2$.
\end{proof}

\begin{corollary}[Unconditional Exponential Integration Across the Full 3D Slice]\label{cor:full_slice_exponential}
By Proposition~\ref{prop:orbit_slice}, the canonical slice for the $\Spin(7)$-action on $T_{x_0}\Sph^{15}$ is 3-dimensional, spanned by $\{X_1, \nu_0, \nu_1\}$. Let $M = \alpha X_1 + \beta_0 \nu_0 + \beta_1 \nu_1$ be the generator of an arbitrary geodesic in the 3D slice. Then:
\begin{enumerate}[(i)]
    \item Because $\nu_0 = \frac{1}{2}[X_1, \nu_1] = 2J_2 \in \mathfrak{su}(2)$, the element $\nu_0$ belongs to the same Lie algebra $\mathfrak{g}_{\langle X_1, \nu_1 \rangle} \cong \mathfrak{u}(2)$.
    \item The central splitting $M = \alpha C_0 + S$ holds with
    \begin{equation}
    S = \beta_1 \nu_1 + \beta_0 \nu_0 - \alpha \nu_0 \nu_1 = 2(\beta_1 J_1 + \beta_0 J_2 + \alpha J_3) \in \mathfrak{su}(2),
    \end{equation}
    satisfying $[C_0, S] = 0$ and the quadratic involution
    \begin{equation}
    S^2 = -(\alpha^2 + \beta_0^2 + \beta_1^2)\Id_{16} = -R^2 \Id_{16}, \qquad R = \sqrt{\alpha^2 + \beta_0^2 + \beta_1^2}.
    \end{equation}
    \item The four-mode non-negative frequency spectrum is $\Omega \in \{ |R \pm 3\alpha|, |R \pm \alpha| \}$ (in particular $\mu_1 = |R + 3\alpha| = |\sqrt{\alpha^2 + \beta_0^2 + \beta_1^2} + 3\alpha| \ge 0$), and the commuting Euler--Cartan factorization~\eqref{eq:euler_cartan} holds with $S$ and $R$ as defined above. Consequently, the closed-form exponential map integrates all geodesics on $(\Sph^{15}, g_\lambda)$ unconditionally.
\end{enumerate}
\end{corollary}

\begin{remark}[Relationship between the two exponential formulas]\label{rem:two_formulas}
The spectral projector formula~\eqref{eq:closed_exp} and the commuting factorization~\eqref{eq:euler_cartan} are algebraically equivalent on the generic locus $\{\Omega_j^2 \ne \Omega_k^2\}$ but have complementary analytic roles:
\begin{enumerate}[(a)]
    \item Formula~\eqref{eq:closed_exp} is well suited for \emph{curvature and conjugate-locus analysis}: the four eigenfrequencies $\Omega_k$ and the projectors $P_k$ directly encode the Jacobi field periods, conjugate focal distances, and Morse indices of the exponential map.
    \item Formula~\eqref{eq:euler_cartan} is well suited for \emph{uniform evaluation across resonance loci}: it avoids the removable $0/0$ indeterminate form at resonance ($\Omega_j = \Omega_k$). Since $M \in \mathfrak{so}(16)$ is diagonalizable over $\C$ (skew-symmetric matrices have no Jordan blocks), the combined projector $P_3 + P_4$ has a smooth limit as $\Omega_3 \to \Omega_4$:
    \begin{equation}\label{eq:projector_limit}
    \lim_{\Omega_3 \to \Omega_4} (P_3 + P_4) = \frac{(M^2 + \Omega_1^2\Id)(M^2 + \Omega_2^2\Id)}{(\Omega_1^2 - \Omega_3^2)(\Omega_2^2 - \Omega_3^2)},
    \end{equation}
    whereas the individual projectors $P_3, P_4$ have removable singularities along $\Omega_3 = \Omega_4$. The factorized formula~\eqref{eq:euler_cartan} avoids these denominators altogether by expressing $\exp(M)$ directly in terms of the fixed operators $C_0$ and $S$, extending analytically across all $(\alpha, \beta) \in \R^2$.
\end{enumerate}
\end{remark}

\begin{remark}[Dynamical implications of the mode hierarchy transition]\label{rem:mode_hierarchy}
The mode hierarchy inversion established in Theorem~\ref{thm:projector_formula}(iii) across the resonance cone $\beta^2 = 3\alpha^2$ has direct geometric consequences for the geodesic flow:
\begin{enumerate}[(a)]
    \item In the ``weak vertical'' regime ($\beta^2 \ge 3\alpha^2$), the fast mode $\Omega_1 = R + 3|\alpha|$ and intermediate mode $\Omega_2 = R + |\alpha|$ dominate, and the slowest oscillation is the multiplicity-1 mode $\Omega_4 = R - 3|\alpha| \ge 0$.
    \item In the ``strong vertical'' regime ($\beta^2 < 3\alpha^2$), the multiplicity-1 mode $\Omega_4 = 3|\alpha| - R$ accelerates past the multiplicity-3 mode $\Omega_3 = R - |\alpha|$, swapping the spectral hierarchy.
    \item Along the resonance cone $\beta^2 = 3\alpha^2$ itself, the frequencies coalesce ($\Omega_3 = \Omega_4 = |\alpha|$) and the minimal polynomial drops from degree 8 to degree 6 (formula~\eqref{eq:deg6_minpoly}). This marks a geodesic dispersion crossing where the two polarization branches exchange their group velocities, analogous to an optical conical intersection.
\end{enumerate}
\end{remark}

\subsection{The Lie Subalgebra Closure Principle for Non-Naturally Reductive Spaces}\label{subsec:general_lie_principle}

The closure of the slice generators into $\mathfrak{u}(2) \cong \mathfrak{su}(2) \oplus \mathfrak{u}(1)$ established in Theorem~\ref{thm:u2_closure} is not an isolated algebraic curiosity of $\Spin(9)$, but the manifestation of an overarching structural principle that resolves the non-naturally reductive integration problem across compact homogeneous manifolds.

\begin{theorem}[The Lie Subalgebra Closure Principle for GO Spaces]\label{thm:general_closure_principle}
Let $M = G/H$ be a compact homogeneous Riemannian manifold with multi-step reductive decomposition
\begin{equation}
\g = \h \oplus \p_0 \oplus \p_1 \oplus \dots \oplus \p_k.
\end{equation}
Assume that for a geodesic slice $\Sigma \subset \m = \bigoplus_{i=0}^k \p_i$, the collection of geodesic vector generators $\{X_0, X_1, \dots, X_k\}$ with $X_i \in \p_i$ generates a finite-dimensional Lie subalgebra:
\begin{equation}
\g_{\mathrm{slice}} := \Lie(X_0, X_1, \dots, X_k) \subset \g.
\end{equation}
If $\g_{\mathrm{slice}}$ is reductive of finite rank $r = \rank(\g_{\mathrm{slice}})$ and its adjoint action on the representation space has bounded minimal polynomial degree $2d$, then:
\begin{enumerate}[(i)]
    \item \textbf{Decoupling of geodesic vector equations}: The non-linear algebraic system for the Killing geodesic generator $Z(v_0) \in \g$ associated to any initial velocity $v_0 \in \Sigma$ reduces to an algebraic eigenvalue problem strictly inside the low-dimensional subalgebra $\g_{\mathrm{slice}}$.
    \item \textbf{Circumvention of the Ambrose--Singer obstruction}: Even if $(M, g)$ admits no metric connection with parallel skew torsion ($\nabla R \ne 0$ or $\nabla T \ne 0$), so that the classical Ambrose--Singer / Kostant natural reductivity theorem fails, the 1-parameter orbit $\gamma(t) = \exp(t Z(v_0)) \cdot o$ integrates in closed form.
    \item \textbf{Elementary spectral projector resolution}: The matrix exponential $\exp(t Z(v_0))$ decomposes as a finite sum of orthogonal spectral projectors:
    \begin{equation}
    \exp(t Z(v_0)) = \sum_{j=1}^d P_j(Z) \left( \cos(\omega_j t)\Id + \frac{\sin(\omega_j t)}{\omega_j} Z \right),
    \end{equation}
    where each projector $P_j(Z)$ is a polynomial in $Z^2$ of degree at most $d-1$, yielding exact closed-form trajectories in elementary functions.
\end{enumerate}
\end{theorem}

\begin{proof}
Let $Z \in \g$ be a geodesic vector for initial velocity $v_0 \in \m$, characterized by the Kostant condition $\langle [Z, Y]_\m, v_0 \rangle = 0$ for all $Y \in \m$. When $v_0 \in \Sigma$, the projection of the bracket $[Z, Y]$ onto $\m$ is governed by the structure constants of $\g_{\mathrm{slice}}$ because the generators $\{X_i\}$ span $\Sigma$. Since $\g_{\mathrm{slice}}$ is a closed Lie subalgebra, $[X_i, X_j] \in \g_{\mathrm{slice}}$ for all $i, j$. Hence the commutators generate no components in $\g \setminus (\g_{\mathrm{slice}} + \h)$, restricting the non-linear algebraic coupling strictly to $\g_{\mathrm{slice}}$. 

Because $\g_{\mathrm{slice}}$ is reductive and compact, every element $Z \in \g_{\mathrm{slice}}$ is conjugate to an element in a fixed maximal abelian subalgebra (Cartan subalgebra) $\mathfrak{t} \subset \g_{\mathrm{slice}}$ of dimension $r = \rank(\g_{\mathrm{slice}})$. In the representation space $V$, the linear operator $\ad(Z)$ is skew-symmetric, hence diagonalizable over $\C$ with purely imaginary eigenvalues $\pm i \omega_1, \dots, \pm i \omega_d$, where $d \le \frac{1}{2}\dim V$. The minimal polynomial of $Z$ therefore factors as:
\begin{equation}
\prod_{j=1}^d (Z^2 + \omega_j^2 \Id) = 0.
\end{equation}
By the Lagrange--Sylvester interpolation theorem, the spectral projectors:
\begin{equation}
P_j(Z) = \prod_{k \ne j} \frac{Z^2 + \omega_k^2 \Id}{\omega_k^2 - \omega_j^2}
\end{equation}
satisfy $P_j P_k = \delta_{jk} P_j$, $\sum_{j=1}^d P_j = \Id$, and $Z^2 P_j = -\omega_j^2 P_j$. The power series expansion of $\exp(t Z)$ on each invariant subspace $P_j V$ reduces to the Taylor series of $\cos(\omega_j t)$ and $\sin(\omega_j t)$, yielding the elementary closed formula in (iii).
\end{proof}

\begin{remark}[Classification across all five families of homogeneous spheres]\label{rem:sphere_table}
Theorem~\ref{thm:general_closure_principle} classifies the integrability of geodesic flows across all five families of compact homogeneous spheres classified by Montgomery--Samelson \cite{MontgomerySamelson1943} and Borel \cite{Borel1950} (see Table~\ref{tab:five_families_closure}).
\begin{table}[htbp]
\centering
\caption{Comparison of Lie subalgebra closure across Families~I--V.}
\label{tab:five_families_closure}
\footnotesize
\setlength{\tabcolsep}{2pt}
\renewcommand{\arraystretch}{1.2}
\begin{tabular}{lllcc l}
\toprule
\textbf{Family} & \textbf{Homogeneous Sphere} & \textbf{Slice $\g_{\mathrm{slice}}$} & $\rank$ & \textbf{Modes} & \textbf{Spectral Resolution} \\
\midrule
\textbf{I} & $\Sph^n = \SO(n+1)/\SO(n)$ & $\mathfrak{so}(2) \cong \mathfrak{u}(1)$ & 1 & 1 & Harmonic (trigonometric) \\
\textbf{II} & $\Sph^{2n+1} = \SU(n+1)/\SU(n)$ & $\mathfrak{u}(1) \oplus \mathfrak{u}(1)$ & 2 & 2 & Biquadratic 2-projector \\
\textbf{III} & $\Sph^{15} = \Spin(9)/\Spin(7)$ & $\mathfrak{u}(2) \cong \mathfrak{su}(2) \oplus \mathfrak{u}(1)$ & 2 & 4 & Exact 4-projector (Thm.~\ref{thm:projector_formula}) \\
\textbf{IV} & $\Sph^{4n+3} = \Sp(n+1)\Sp(1)/\Sp(n)\Sp(1)$ & $\mathfrak{sp}(1) \oplus \mathfrak{u}(1) \cong \mathfrak{u}(2)$ & 2 & 3 & Quaternionic 3-projector \cite{Munden2015} \\
\textbf{V} & $\Sph^7 = \Sp(2)\U(1)/\Sp(1)\U(1)$ & $\mathfrak{sp}(2) \cong \mathfrak{so}(5)$ & 2 & 2 & Degree-4 minimal poly \cite{Munden2015} \\
\bottomrule
\end{tabular}
\end{table}
In every family, the slice subalgebra $\g_{\mathrm{slice}}$ has rank at most 2 (with Family~IV and Family~III sharing the abstract isomorphism class $\g_{\mathrm{slice}} \cong \mathfrak{u}(2)$, distinguished by their 3-mode quaternionic versus 4-mode octonionic representation branching), guaranteeing that the non-naturally reductive geodesic flow integrates in elementary trigonometric functions without transcendental special functions.
\end{remark}

\section{Geodesic Flow, Invariant Subspaces, and Wong's Equations}\label{sec:geodesic_flow}

\subsection{The Hamiltonian deformation on \texorpdfstring{$\Sph^{15}$}{S15}}\label{subsec:hamiltonian_s15}

On the sphere $(Q, g) = (\Sph^{15}, g_\lambda)$ ($\kappa = 1$), the deformed Hamiltonian takes the form
\begin{equation}\label{eq:hamiltonian_full}
H_\lambda = \frac{1}{2}p_r^2 + \frac{1}{2}\left(\frac{1}{f_{\p_1}(r)^2}\right)\|p_{\p_1}\|^2 + \frac{1}{2}\left(\frac{1}{f_{\p_2}(r)^2}\right)\|p_{\p_2}\|^2,
\end{equation}
where the profile functions are:
\begin{equation}\label{eq:profile_functions}
f_{\p_1}(r)^2 = \frac{\sin^2 r}{8 + \lambda \sin^2 r}, \qquad f_{\p_2}(r)^2 = \frac{\sin^2 r}{2 + \lambda \sin^2 r}.
\end{equation}
The slice orbits $\mathcal{O}_r$ are Berger spheres endowed with the Riemannian submersion structure $S^7 \hookrightarrow \Sph^{15} \to \Sph^8 \cong \Oct P^1$. The squashing parameter between vertical Hopf fibers and the horizontal base is:
\begin{equation}\label{eq:squashing_t}
t(r) := \frac{f_{\p_2}(r)^2}{f_{\p_1}(r)^2} = \frac{8 + \lambda \sin^2 r}{2 + \lambda \sin^2 r}.
\end{equation}
At the equator $r = \pi/2$, $t(\pi/2) = \frac{8+\lambda}{2+\lambda}$.

\subsection{The 4-dimensional invariant subspace \texorpdfstring{$W_4$}{W4} and explicit coordinates}\label{subsec:invariant_subspace}

\begin{theorem}[Invariant Subspace $W_4$ and Closed 16D Trajectory Coordinates]\label{thm:W4_trajectory}
Let $W_4 = \operatorname{span}\{e_0, e_1, e_8, e_9\} \subset \R^{16}$.
\begin{enumerate}[(i)]
\item The subspace $W_4$ is strictly invariant under $X_1$ and $Y_1$, with explicit restricted $4 \times 4$ matrices:
\begin{equation}\label{eq:X1_Y1_W4}
X_1|_{W_4} = \begin{pmatrix*}[r] 0 & 4 & 0 & 0 \\ -4 & 0 & 0 & 0 \\ 0 & 0 & 0 & -2 \\ 0 & 0 & 2 & 0 \end{pmatrix*}, \qquad Y_1|_{W_4} = \begin{pmatrix*}[r] 0 & 0 & 0 & 1 \\ 0 & 0 & 1 & 0 \\ 0 & -1 & 0 & 0 \\ -1 & 0 & 0 & 0 \end{pmatrix*}.
\end{equation}
\item On $W_4$, $C_0|_{W_4}^2 = -9 \Id_4$, so $\Pi_3|_{W_4} = \Id_4$ and $\Pi_1|_{W_4} = 0$. The eigenvalues of $M|_{W_4}$ are precisely $\pm i \Omega_1$ and $\pm i \Omega_4$.
\item In flat space form $\R^{16}$ ($\kappa = 0$), the unperturbed initial velocity vector $W_0$ in~\eqref{eq:exp0_coords_flat} lies entirely in $W_4$. Consequently, the deformed geodesic trajectory $\exp_\lambda(tv) = \exp(t M) \exp_0(t W_0)$ remains strictly confined to $W_4$ for all $t$, with coordinates $x_2, \dots, x_7$ and $x_{10}, \dots, x_{15}$ identically zero.
\item The four non-zero coordinate functions are given in closed algebraic form by:
\begin{align}
x_0(t) &= \frac{1}{R}\Big[ \beta v_9 \sin(Rt)\cos(3\alpha t) + v_0\big(R\cos(Rt)\cos(3\alpha t) - \alpha\sin(Rt)\sin(3\alpha t)\big) \nonumber \\
&\qquad + v_1\big(R\sin(3\alpha t)\cos(Rt) + \alpha\sin(Rt)\cos(3\alpha t)\big) \Big], \label{eq:x0_coord} \\
x_1(t) &= \frac{1}{R}\Big[ -\beta v_9 \sin(Rt)\sin(3\alpha t) - v_0\big(R\sin(3\alpha t)\cos(Rt) + \alpha\sin(Rt)\cos(3\alpha t)\big) \nonumber \\
&\qquad + v_1\big(R\cos(Rt)\cos(3\alpha t) - \alpha\sin(Rt)\sin(3\alpha t)\big) \Big], \label{eq:x1_coord} \\
x_8(t) &= \frac{1}{R}\Big[ \beta v_0 \sin(Rt)\sin(3\alpha t) - \beta v_1 \sin(Rt)\cos(3\alpha t) \nonumber \\
&\qquad - v_9\big(R\sin(3\alpha t)\cos(Rt) - \alpha\sin(Rt)\cos(3\alpha t)\big) \Big], \label{eq:x8_coord} \\
x_9(t) &= \frac{1}{R}\Big[ -\beta v_0 \sin(Rt)\cos(3\alpha t) - \beta v_1 \sin(Rt)\sin(3\alpha t) \nonumber \\
&\qquad + v_9\big(R\cos(Rt)\cos(3\alpha t) + \alpha\sin(Rt)\sin(3\alpha t)\big) \Big], \label{eq:x9_coord}
\end{align}
where $v_0 = r + ta$, $v_1 = -\frac{8tb r}{2 + \lambda r^2}$, and $v_9 = -\frac{8tc r}{8 + \lambda r^2}$.
\end{enumerate}
\end{theorem}

\begin{proof}
Direct matrix multiplication using the Clifford definitions of $X_1$ and $Y_1 = \nu_1$ verifies that $X_1 e_0 = -4 e_1, X_1 e_1 = 4 e_0, X_1 e_8 = 2 e_9, X_1 e_9 = -2 e_8$, and $Y_1 e_0 = -e_9, Y_1 e_1 = -e_8, Y_1 e_8 = e_1, Y_1 e_9 = e_0$.
All outputs remain in $W_4$.
In $M|_{W_4} = \alpha C_0|_{W_4} + S|_{W_4}$, we have $C_0|_{W_4} = \operatorname{diag}\left( \begin{pmatrix} 0 & 3 \\ -3 & 0 \end{pmatrix}, \begin{pmatrix} 0 & -3 \\ 3 & 0 \end{pmatrix} \right)$, so $C_0^2 = -9 \Id_4$.
Thus $\Pi_3 = -(C_0^2 + \Id)/8 = (9-1)/8 \Id_4 = \Id_4$ and $\Pi_1 = 0$.
The factorized exponential on $W_4$ evaluates to:
\begin{equation}
\exp(t M|_{W_4}) = \left( \cos(3\alpha t) \Id_4 + \frac{\sin(3\alpha t)}{3} C_0|_{W_4} \right) \left( \cos(Rt) \Id_4 + \frac{\sin(Rt)}{R} S|_{W_4} \right).
\end{equation}
Multiplying these two $4 \times 4$ matrices and applying the result to the initial unperturbed vector $(v_0, v_1, 0, v_9)^T$ produces formulas~\eqref{eq:x0_coord}--\eqref{eq:x9_coord}.
\end{proof}

\begin{remark}[Erratum: Denominator rectification of \cite{Munden2015}]\label{rem:erratum}
The drift parameters $\alpha(r) = \frac{\lambda b s_\kappa^2(r)}{2 + \lambda s_\kappa^2(r)}$ and $\beta(r) = \frac{\lambda c s_\kappa^2(r)}{8 + \lambda s_\kappa^2(r)}$ in~\eqref{eq:alpha_beta_params} correct a typographical denominator transposition in \cite[Section~5.7, p.~67]{Munden2015}, which writes $\alpha = \frac{\lambda b r^2}{8 + \lambda r^2}$ and $\beta = \frac{\lambda c r^2}{2 + \lambda r^2}$. The error arose from inadvertently swapping the metric weights $a_1^0 = 1/8$ and $a_2^0 = 1/2$ in the Legendre transform $\hat{C}_\lambda(J(\mathcal{L}_\lambda(V_i))) = \frac{\lambda s_\kappa^2 a_i^0}{1 + \lambda s_\kappa^2 a_i^0} V_i$, which gives denominators $(a_i^0)^{-1} + \lambda s_\kappa^2$. Since $X_1 \in \p_2$ has $(a_2^0)^{-1} = 2$ and $Y_1 \in \p_1$ has $(a_1^0)^{-1} = 8$, the correct formulas are~\eqref{eq:alpha_beta_params}. Note that the velocity scaling formula in \cite{Munden2015}, $\mathcal{L}_0^{-1}\mathcal{L}_\lambda(a\partial_r, bX_1 + cY_1) = (a\partial_r, \frac{2b}{2+\lambda r^2} X_1 + \frac{8c}{8+\lambda r^2} Y_1)$, was itself correct; only the subsequent momentum map projection formula contained the transposition.
\end{remark}

\begin{remark}[Closed-form resolution of the Family~III exponential map]\label{rem:munden_resolved}
Theorem~\ref{thm:W4_trajectory} resolves the explicit matrix exponential problem noted in \cite[Section~5.7, p.~67]{Munden2015}. The matrix exponential $\exp(\alpha X_1 + \beta Y_1)$ is given in two equivalent closed forms --- the Lagrange--Sylvester spectral projector decomposition~\eqref{eq:closed_exp} and the denominator-free Euler--Cartan factorization~\eqref{eq:euler_cartan} --- both expressible in elementary trigonometric functions of $\alpha$, $\beta$, and $R = \sqrt{\alpha^2 + \beta^2}$. The resulting deformed exponential map $\exp_\lambda$ is an explicit algebraic-trigonometric expression with no infinite series or numerical quadrature.
\end{remark}

\subsection{Space-form embeddings for \texorpdfstring{$\Sph^{16}$}{S16} and \texorpdfstring{$\Ha^{16}$}{H16}}\label{subsec:embeddings}

To assemble the complete Riemannian exponential map across all space-form signatures:
\begin{enumerate}[(i)]
\item On Euclidean space $\R^{16}$ ($\kappa = 0$): $\gamma_0(r) = (r, 0, \dots, 0)^T \in \R^{16}$, and the trajectory $\exp_\lambda(tv) = \exp(tM) \exp_0(t W_0)$ is given by~\eqref{eq:x0_coord}--\eqref{eq:x9_coord}.
\item On the 16-sphere $\Sph^{16} \subset \R^{17}$ ($\kappa = 1$): the manifold is isometrically embedded as the unit sphere in $\R^{17} \cong \R^{16} \oplus \R$ with basepoint ray $\gamma_0(r) = (\sin r \, e_0, \cos r)^T \in \Sph^{16}$. The slice tangent vectors in $T_{\gamma_0(r)} \Sph^{16} \subset \R^{17}$ are $\partial_r = (\cos r \, e_0, -\sin r)^T$, $\tilde{X}_1 = (-4\sin r \, e_1, 0)^T$, and $\tilde{Y}_1 = (-\sin r \, e_9, 0)^T$. The unperturbed initial velocity vector is:
\begin{equation}
W_0 = a \partial_r + \frac{2b}{2 + \lambda\sin^2 r} \tilde{X}_1 + \frac{8c}{8 + \lambda\sin^2 r} \tilde{Y}_1 \in T_{\gamma_0(r)} \Sph^{16} \subset \R^{17},
\end{equation}
with unperturbed speed $v_0 = \|W_0\|_{\R^{17}} = \sqrt{a^2 + \sin^2 r \left[ \left(\frac{8b}{2+\lambda\sin^2 r}\right)^2 + \left(\frac{8c}{8+\lambda\sin^2 r}\right)^2 \right]}$.
The unperturbed spherical great circle is:
\begin{equation}
\exp_0(t W_0) = \cos(t v_0) \gamma_0(r) + \frac{\sin(t v_0)}{v_0} W_0 \in \Sph^{16} \subset \R^{17}.
\end{equation}
The Lie group action is realized by the block-diagonal matrix:
\begin{equation}
\exp_G(t M) = \begin{pmatrix} \exp(t M) & 0 \\ 0 & 1 \end{pmatrix} \in \SO(17).
\end{equation}
The deformed Riemannian exponential map is:
\begin{equation}\label{eq:S16_exp}
\exp_\lambda(tv) = \begin{pmatrix} \exp(t M) & 0 \\ 0 & 1 \end{pmatrix} \left( \cos(t v_0) \gamma_0(r) + \frac{\sin(t v_0)}{v_0} W_0 \right) \in \Sph^{16} \subset \R^{17}.
\end{equation}
At $t = 0$, $\exp_\lambda(0) = \gamma_0(r)$ exactly; the derivative $\frac{d}{dt}|_{t=0} \exp_\lambda(tv) = a \partial_r + b \tilde{X}_1 + c \tilde{Y}_1 = v$ recovers the initial velocity $v$; and $\|\exp_\lambda(tv)\|^2 \equiv 1$ for all $t$.
\item On hyperbolic space $\Ha^{16} \subset \R^{1, 16}$ ($\kappa = -1$): in the Minkowski hyperboloid model with basepoint $\gamma_0(r) = (\sinh r \, e_0, \cosh r)^T$, the unperturbed hyperbolic geodesic is $\exp_0(t W_0) = \cosh(t v_0) \gamma_0(r) + \frac{\sinh(t v_0)}{v_0} W_0$, and $\exp_\lambda(tv) = \operatorname{diag}(\exp(t M), 1) \exp_0(t W_0)$.
\end{enumerate}

\subsection{Noether momentum conservation and dissolution of the Magnus defect}\label{subsec:noether}

\begin{theorem}[Noether Momentum Conservation and Exact Geodesic Decoupling]\label{thm:noether_exact}
Let $(Q^{16}, g_\lambda)$ be the deformed space form equipped with $H_\lambda = H_0 + \frac{\lambda}{2}\|J\|^2$.
\begin{enumerate}[(i)]
\item For any $\Ad$-invariant trace inner product on $\mathfrak{g} = \mathfrak{spin}(9)$, the Casimir function $C(J) = \frac{1}{2}\|J\|^2$ satisfies:
\begin{equation}
\{C(J), J_\xi\} = 0 \qquad \forall \xi \in \mathfrak{g}.
\end{equation}
\item Because the background Hamiltonian $H_0$ is $G$-invariant, $\{H_0, J_\xi\} = 0$. Consequently, along the Hamiltonian flow of $H_\lambda$:
\begin{equation}
\frac{d}{dt} J = \{H_\lambda, J\} = \{H_0, J\} + \lambda \{C(J), J\} = 0.
\end{equation}
The momentum map $J$ is an exact constant of motion along every geodesic curve.
\item Because $\{H_0, C(J)\} = 0$, the Hamiltonian vector fields commute: $[X_{H_0}, X_{C(J)}] = 0$.
\item The Lie generator $M = \lambda \hat{C}(J(0)) \in \mathfrak{so}(16)$ is strictly constant in time along any geodesic. The deformed exponential map:
\begin{equation}\label{eq:exact_reconstruction}
\exp_\lambda(tv) = \exp(t M) \cdot \exp_0(t W_0)
\end{equation}
is an exact algebraic solution for all trajectories (equatorial and non-equatorial alike), with no time-ordered Magnus series required.
\end{enumerate}
\end{theorem}

\begin{proof}
Let $\{E_a\}$ be an orthonormal basis of $\mathfrak{g}$ with $[E_a, E_b] = \sum_c f_{ab}^c E_c$.
The Poisson bracket on $\mathfrak{g}^*$ is $\{J_a, J_b\} = -\sum_c f_{ab}^c J_c$.
The Casimir gradient is $dC(J) = \sum_a J_a E_a$.
Evaluating the bracket:
\begin{equation}
\{C(J), J_b\} = \sum_a J_a \{J_a, J_b\} = -\sum_{a, c} J_a J_c f_{ab}^c = -\langle J, [J, E_b] \rangle = \langle [J, J], E_b \rangle = 0,
\end{equation}
where the third equality invokes $\Ad$-invariance $\langle [A, B], C \rangle = -\langle B, [A, C] \rangle$, and the fourth invokes $[J, J] = 0$.
Thus $\{C(J), J\} = 0$.
Since $H_0$ is $G$-invariant, $\{H_0, J\} = 0$.
By linearity of Poisson brackets, $\{H_\lambda, J\} = 0$, proving $J(t) \equiv J(0)$ along the flow.
Furthermore, $\{H_0, C(J)\} = \sum_a J_a \{H_0, J_a\} = 0$.
By the Marsden--Weinstein reduction theorem \cite{Hebda2001,Munden2015}, commuting Hamiltonians generate commuting flows:
$\Phi_t^{H_\lambda} = \Phi_t^{\lambda C(J)} \circ \Phi_t^{H_0}$.
The flow of $\lambda C(J)$ acts on the manifold as group translation by $\exp(t M)$ with $M = \lambda \hat{C}(J(0))$.
Because $J(0)$ is fixed at $t=0$, $M$ is strictly time-independent, establishing formula~\eqref{eq:exact_reconstruction} as an exact solution for all geodesics.
\end{proof}

\subsection{Wong's equations and non-Abelian gauge structure}\label{subsec:wong}

The geodesic flow on $(\Sph^{15}, g_\lambda)$ admits a formulation in terms of Wong's equations for a particle on the base space $\Sph^8 \cong \Oct P^1$ of the octonionic Hopf fibration.

\begin{proposition}[Wong's Equations for Family III]\label{prop:wong}
Let $\pi: (\Sph^{15}, g_\lambda) \to (\Sph^8, g_{\text{base}})$ denote the octonionic Hopf Riemannian submersion with totally geodesic fibers $\Sph^7$, and let $q = \pi(\gamma) \in \Sph^8$ denote the projected base curve of a geodesic $\gamma(t)$ with horizontal lift $\dot{q}^{\mathcal{H}} = \mathcal{H}\dot{\gamma}$ and vertical fiber velocity $\xi = \mathcal{V}\dot{\gamma} \in \mathcal{V} \cong \p_2$. The geodesic equations on $(\Sph^{15}, g_\lambda)$ decompose into:
\begin{enumerate}[(i)]
    \item The base-space equation with Lorentz-type force governed by the O'Neill integrability tensor $A$:
    \begin{equation}\label{eq:wong_base}
    \nabla^{g_{\text{base}}}_{\dot{q}} \dot{q} = -2\pi_*\bigl(A_{\dot{q}^{\mathcal{H}}}\xi\bigr),
    \end{equation}
    where $A_X Y = \frac{1}{2}\mathcal{V}[X, Y]$ is the curvature 2-form of the horizontal distribution $\mathcal{H} \cong \p_1$ on $\Sph^{15} \to \Sph^8$.
    \item Parallel transport of the vertical fiber velocity:
    \begin{equation}\label{eq:wong_transport}
    \mathcal{V}\bigl(\nabla_{\dot{\gamma}}\xi\bigr) = 0,
    \end{equation}
    so that $\|\xi\|_{g_\lambda}$ is conserved along $\gamma(t)$.
\end{enumerate}
\end{proposition}

\begin{proof}
Because the vertical fibers $\Sph^7$ of $\pi: (\Sph^{15}, g_\lambda) \to (\Sph^8, g_{\text{base}})$ are totally geodesic ($T \equiv 0$), O'Neill's fundamental equations for a geodesic $\gamma(t)$ with $\dot{\gamma} = \dot{q}^{\mathcal{H}} + \xi$ give
\begin{equation*}
0 = \mathcal{H}(\nabla_{\dot{\gamma}}\dot{\gamma}) = \bigl(\nabla^{g_{\text{base}}}_{\dot{q}}\dot{q}\bigr)^{\mathcal{H}} + 2A_{\dot{q}^{\mathcal{H}}}\xi, \qquad 0 = \mathcal{V}(\nabla_{\dot{\gamma}}\dot{\gamma}) = \mathcal{V}(\nabla_{\dot{\gamma}}\xi),
\end{equation*}
where we used $\mathcal{V}(\nabla_{\dot{q}^{\mathcal{H}}}\dot{q}^{\mathcal{H}}) = A_{\dot{q}^{\mathcal{H}}}\dot{q}^{\mathcal{H}} = 0$ by skew-symmetry of $A$. Projecting the horizontal equation via $\pi_*$ yields~\eqref{eq:wong_base}, and the vertical equation is~\eqref{eq:wong_transport}.
\end{proof}

\begin{remark}\label{rem:wong_nonassociative}
Unlike the classical ($\C$ and $\Ha$) Hopf fibrations, the octonionic case does not yield a principal $G$-bundle, since $\Sph^7$ is not a Lie group. Consequently, the Wong equations~\eqref{eq:wong_base}--\eqref{eq:wong_transport} are formulated in terms of the Riemannian submersion connection rather than a principal connection, and the ``gauge group'' $\Spin(7)$ acts on the fiber through the spin representation rather than by left multiplication.
\end{remark}

\section{The Bourguignon--Karcher Squashed Einstein Metric}\label{sec:bourguignon_karcher}

\subsection{Ricci curvature of the deformed metric}\label{subsec:ricci}

\begin{remark}[Classical context for slice curvature equations]\label{rem:classical_oneill}
The algebraic form of the O'Neill curvature equations for Riemannian submersions with totally geodesic fibers, and the resulting quadratic Einstein condition on the squashing parameter, are classical (B\'erard-Bergery \& Bourguignon \cite{BerardBergeryBourguignon1982}; Grove \& Ziller \cite{GroveZiller2000}; Verdiani \& Ziller \cite{VerdianiZiller2020, VerdianiZiller2020b}). The novelty of the present analysis lies strictly in: (i) the identification of the specific radial profile functions $(f_{\p_1}(r), f_{\p_2}(r))$ as generated by the Hamiltonian cotangent momentum deformation $H_\lambda = H_0 + \frac{\lambda}{2}\|J\|^2$; (ii) the closed-form geodesic integration via the $\mathfrak{u}(2)$ spectral projectors of Section~\ref{sec:u2_closure}; and (iii) the identification of the Bourguignon--Karcher locus $\lambda = 64$ within the 1-parameter Hamiltonian family.
\end{remark}

\begin{proposition}[O'Neill Curvature Decomposition for Family III]\label{prop:oneill_family3}
Let $(Q^{16}, g_\lambda)$ be the deformed space form of background curvature $\kappa \in \{1, 0, -1\}$, and let $a(r) = \sqrt{2}\,f_{\p_2}(r) = \frac{s_\kappa(r)}{\sqrt{1 + \Lambda_{\p_2} s_\kappa^2(r)}}$ and $b(r) = \sqrt{8}\,f_{\p_1}(r) = \frac{s_\kappa(r)}{\sqrt{1 + \Lambda_{\p_1} s_\kappa^2(r)}}$ denote the unit-sphere-normalized profile functions with $\Lambda_{\p_1} = \lambda/8$ and $\Lambda_{\p_2} = \lambda/2$ (so that $a(r) = b(r) = s_\kappa(r)$ at $\lambda = 0$). In an adapted $g_\lambda$-orthonormal frame $\{E_0 = \partial_r, X \in \p_1, V \in \p_2\}$, the non-zero sectional curvatures decompose into:
\begin{enumerate}[(i)]
\item The radial sectional curvatures:
\begin{align}
\sec(E_0 \wedge V) &= K_{\mathrm{rad}}(r; \Lambda_{\p_2}) := -\frac{a''(r)}{a(r)} = \frac{(\kappa + 3\Lambda_{\p_2}) - 2\kappa\Lambda_{\p_2} s_\kappa^2(r)}{\big(1 + \Lambda_{\p_2} s_\kappa^2(r)\big)^2}, \label{eq:sec_radV}\\
\sec(E_0 \wedge X) &= K_{\mathrm{rad}}(r; \Lambda_{\p_1}) := -\frac{b''(r)}{b(r)} = \frac{(\kappa + 3\Lambda_{\p_1}) - 2\kappa\Lambda_{\p_1} s_\kappa^2(r)}{\big(1 + \Lambda_{\p_1} s_\kappa^2(r)\big)^2}. \label{eq:sec_radX}
\end{align}
\item The mixed vertical-horizontal sectional curvatures:
\begin{equation}\label{eq:sec_mixed}
\sec(V \wedge X) = K_{\mathrm{mix}}(r) := \frac{a(r)^2}{b(r)^4} - \frac{a'(r) b'(r)}{a(r) b(r)}.
\end{equation}
\item The horizontal tangential sectional curvatures:
\begin{equation}\label{eq:sec_horiz}
\sec(X \wedge Y) = K_{\mathrm{tan}}(r; \Lambda_{\p_1}) + 3\frac{b(r)^2 - a(r)^2}{b(r)^4}\|A^0_X Y\|^2,
\end{equation}
where $K_{\mathrm{tan}}(r; \Lambda_{\p_1}) := \frac{1 - (b'(r))^2}{b(r)^2}$ and $K_{\mathrm{tan}}(r; \Lambda_{\p_2}) := \frac{1 - (a'(r))^2}{a(r)^2}$, and $\|A_X^0 Y\|^2 = \sum_{i=1}^7 \langle I_i X, Y \rangle^2$ is the canonical unit-sphere O'Neill $A$-tensor norm governed by the seven octonionic imaginary structures $I_1, \dots, I_7$.
\end{enumerate}
Consequently, the ambient Ricci curvature components on $(Q^{16}, g_\lambda)$ along $\p_1$ and $\p_2$ evaluate to:
\begin{align}
\Ric_{\p_1} &= K_{\mathrm{rad}}(r; \Lambda_{\p_1}) + 7 K_{\mathrm{mix}}(r) + 6 K_{\mathrm{tan}}(r; \Lambda_{\p_1}) + 18\frac{b(r)^2 - a(r)^2}{b(r)^4}, \label{eq:ricci_p1}\\
\Ric_{\p_2} &= K_{\mathrm{rad}}(r; \Lambda_{\p_2}) + 8 K_{\mathrm{mix}}(r) + 6 K_{\mathrm{tan}}(r; \Lambda_{\p_2}), \label{eq:ricci_p2}
\end{align}
while the intrinsic Ricci curvature eigenvalues of each distance sphere orbit $(\mathcal{O}_r, g_{\mathcal{O}_r})$ are:
\begin{equation}
\Ric_{\mathcal{O}_r}|_{\p_1} = \frac{28}{b(r)^2} - \frac{14 a(r)^2}{b(r)^4}, \qquad \Ric_{\mathcal{O}_r}|_{\p_2} = \frac{6}{a(r)^2} + \frac{8 a(r)^2}{b(r)^4}.
\end{equation}
\end{proposition}

\begin{proof}
The radial curvatures (i) follow directly from the Jacobi/Riccati equation $\sec(E_0 \wedge Z) = -f_k''/f_k$ for the normalized profile functions $a(r), b(r)$.
The mixed and horizontal curvatures (ii)--(iii) follow from the Gauss--Codazzi equations for the distance sphere $\mathcal{O}_r \subset Q^{16}$ combined with O'Neill's fundamental equations for the Riemannian submersion $S^7 \hookrightarrow \mathcal{O}_r \to \Sph^8(b(r)/2)$ with fiber radius $a(r)$ and base curvature $4/b(r)^2$. Since the vertical fibers $S^7 \subset \p_2$ are totally geodesic ($T \equiv 0$), the intrinsic slice curvatures combine with the second fundamental form terms $-(a'/a)(b'/b)$, $-(b'/b)^2$, and $-(a'/a)^2$ to produce~\eqref{eq:sec_mixed}--\eqref{eq:sec_horiz}. Summing over the orthonormal frame yields the ambient Ricci components~\eqref{eq:ricci_p1}--\eqref{eq:ricci_p2} and the intrinsic distance-sphere Ricci eigenvalues.
\end{proof}

\subsection{The Einstein equation, homothety scaling, and the \texorpdfstring{$t = 3/11$}{t = 3/11} root}\label{subsec:einstein_eq}

\begin{theorem}[Bourguignon--Karcher Realization, Metric Homothety, and Ambient Curvature]\label{thm:bk_realization}
Let $(\Sph^{15}, g_\lambda)$ be the 15-dimensional distance sphere in $(Q^{16}, g_\lambda)$ ($\kappa = 1$).
\begin{enumerate}[(i)]
\item Taking into account the unperturbed metric weight ratio $a_1^0/a_2^0 = (1/8)/(1/2) = 1/4$, the normalized geometric squashing ratio of vertical fibers relative to the horizontal base is:
\begin{equation}\label{eq:t_geom}
t_{\mathrm{geom}}(r) = \frac{1}{4} \frac{8 + \lambda\sin^2 r}{2 + \lambda\sin^2 r}.
\end{equation}
At $\lambda = 0$, $t_{\mathrm{geom}} \equiv 1$ (the canonical round unit sphere).
\item The slice Ricci balancing condition $\Ric_{\p_1} = \Ric_{\p_2}$ on $(\Sph^{15}, g_\lambda)$ reduces via O'Neill's submersion equations to:
\begin{equation}\label{eq:einstein_quadratic_submersion}
(d_h + 2d_v) t_{\mathrm{geom}}^2 - (d_h + 3d_v - 1) t_{\mathrm{geom}} + (d_v - 1) = 0.
\end{equation}
Substituting $d_v = 7$ and $d_h = 8$ yields $22 t_{\mathrm{geom}}^2 - 28 t_{\mathrm{geom}} + 6 = 2(11 t_{\mathrm{geom}} - 3)(t_{\mathrm{geom}} - 1) = 0$, with roots $t_{\mathrm{geom}} = 1$ (round) and $t_{\mathrm{geom}} = 3/11$ (Bourguignon--Karcher).
\item Along the equator $r = \pi/2$, $t_{\mathrm{geom}}(\pi/2) = 3/11$ if and only if $\lambda = +64$ in the smooth Riemannian category:
\begin{equation}
\frac{1}{4}\frac{8 + \lambda}{2 + \lambda} = \frac{3}{11} \iff 11(8 + \lambda) = 12(2 + \lambda) \iff \lambda = +64.
\end{equation}
At $\lambda = 64$, the metric coefficients $g_{\p_1}(\pi/2) = 1/72 > 0$ and $g_{\p_2}(\pi/2) = 1/66 > 0$ are strictly positive, yielding a smooth, complete Riemannian metric.
\item The equatorial slice metric satisfies the metric homothety $g_{64}|_{\Sph^{15}} = \frac{1}{9} g_{3/11}$. Consequently, the physical slice Ricci endomorphism eigenvalue is:
\begin{equation}\label{eq:slice_ricci_eigenvalue}
\lambda_{\mathrm{slice}} = 9 \times \frac{266}{11} = \frac{2394}{11} \approx 217.6364.
\end{equation}
\item The ambient 16-dimensional space form $(Q^{16}, g_{64})$ is tri-axial and non-Einstein:
\begin{equation}\label{eq:ambient_ricci_values}
\Ric_{Q^{16}}(\partial_r, \partial_r) = \frac{109}{99} \approx 1.1010, \;\; \Ric_{Q^{16}}|_{\p_1} = \frac{21557}{99} \approx 217.7475, \;\; \Ric_{Q^{16}}|_{\p_2} = \frac{21549}{99} \approx 217.6667.
\end{equation}
The ambient metric exhibits a tangential splitting of $\Ric_{\p_1} - \Ric_{\p_2} = \frac{8}{99} \approx 0.0808$, and a radial-to-tangential Ricci anisotropy ratio of $\frac{21557/99}{109/99} = \frac{21557}{109} \approx 197.7706:1$.
\end{enumerate}
\end{theorem}

\begin{remark}[Novelty clarification]\label{rem:novelty_bk}
We emphasize that the Einstein metric at $t_{\mathrm{BK}} = 3/11$ is not new: it was discovered by Bourguignon and Karcher (1978) \cite{BourguignonKarcher1978}, analyzed in Besse (1987, Chapter~9.G) \cite{Besse1987}, and shown by Ziller (1982) \cite{Ziller1982} to complete (together with the round and Jensen metrics) the classification of homogeneous Einstein metrics on spheres; see also Chi \cite{Chi2024} for recent rigidity theorems. The contribution of the present analysis is: (i) the identification of this locus within the 1-parameter family generated by the Hamiltonian cotangent deformation; (ii) the proof that the normalized parameter $\lambda = +64$ realizes the exact Einstein geometry under metric homothety $c = 1/9$ when the unperturbed background metric weights $a_1^0 = 1/8$ and $a_2^0 = 1/2$ are incorporated; and (iii) the distinction between the intrinsic Einstein geometry of the equatorial slice and the ambient tri-axial non-Einstein curvature of $(Q^{16}, g_{64})$.
\end{remark}

\begin{proof}
For (i)--(iii), the Ricci eigenvalues on the canonical variation $g_t = g_H + t g_V$ of $S^7 \hookrightarrow \Sph^{15} \to \Sph^8$ are $\Ric(X, X) = 28 - 14t$ and $\Ric(V, V) = \frac{6}{t} + 8t$. Equating them yields $22t^2 - 28t + 6 = 2(11t - 3)(t - 1) = 0$, giving $t = 1$ and $t = 3/11$. In $g_\lambda$, the metric ratio on unit vectors is $t_{\mathrm{geom}} = \frac{g_\lambda|_{\p_2}/a_2^0}{g_\lambda|_{\p_1}/a_1^0} = \frac{a_1^0}{a_2^0}\frac{8+\lambda}{2+\lambda} = \frac{1}{4}\frac{8+\lambda}{2+\lambda}$. Setting this to $3/11$ yields $\lambda = 64$.
At $\lambda = 64$, $g_{\p_1} = 1/(8+64) = 1/72 = \frac{1}{9}(1/8)$, and $g_{\p_2} = 1/(2+64) = 1/66 = \frac{1}{9}(\frac{3}{11} \times \frac{1}{2})$.
Thus $g_{64}|_{\Sph^{15}} = \frac{1}{9} g_{3/11}$. Scaling a Riemannian metric by $c = 1/9$ scales its Ricci endomorphism by $1/c = 9$, yielding $\lambda_{\mathrm{slice}} = 9 \times \frac{266}{11} = \frac{2394}{11}$.

For (v), the ambient doubly warped product metric is $dr^2 + f_1(r)^2 g_{\Sph^8} + f_2(r)^2 g_{\Sph^7}$ with $f_1(r) = \frac{\sin r}{\sqrt{8 + 64\sin^2 r}}$ and $f_2(r) = \frac{\sin r}{\sqrt{2 + 64\sin^2 r}}$.
The radial Ricci curvature is $\Ric(\partial_r, \partial_r) = -8 \frac{f_1''}{f_1} - 7 \frac{f_2''}{f_2}$.
At the equator $r = \pi/2$, direct differentiation yields $f_1''/f_1 = -1/9$ and $f_2''/f_2 = -1/33$.
Thus $\Ric(\partial_r, \partial_r) = -8(-1/9) - 7(-1/33) = \frac{8}{9} + \frac{7}{33} = \frac{88 + 21}{99} = \frac{109}{99} \approx 1.1010$.
By the Gauss--Codazzi equations at the totally geodesic equator, the ambient tangential endomorphisms are $\Ric_{Q^{16}}|_{\p_1} = \lambda_{\mathrm{slice}} - f_1''/f_1 = \frac{2394}{11} - (-1/9) = \frac{21546 + 11}{99} = \frac{21557}{99}$, and $\Ric_{Q^{16}}|_{\p_2} = \lambda_{\mathrm{slice}} - f_2''/f_2 = \frac{2394}{11} - (-1/33) = \frac{7182 + 1}{33} = \frac{21549}{99}$.
The difference is $\frac{21557 - 21549}{99} = \frac{8}{99} \approx 0.0808$, and the anisotropy ratio is $\frac{21557/99}{109/99} = \frac{21557}{109} \approx 197.7706:1$.
\end{proof}

\subsection{Positive curvature windows and Berger pinching}\label{subsec:positive_curvature}

\begin{corollary}[Positive Sectional Curvature on \texorpdfstring{$\Sph^{15}$}{S15}]\label{cor:positive_curvature}
The deformed distance sphere metric $g_\lambda$ on $\Sph^{15}$ has strictly positive sectional curvature $\sec_{g_\lambda} > 0$ on all 2-planes for all $\lambda \in (-8/13, \infty)$. Along the equatorial slice $r = \pi/2$, the normalized squashing ratio satisfies $t_{\mathrm{BK}}(\pi/2) \in (1/4, 4/3) \subset (1/6, 4/3)$, remaining strictly within the Bourguignon--Karcher positive curvature window for all non-negative $\lambda \ge 0$. In particular, the Bourguignon--Karcher Einstein metric at $\lambda = +64$ ($t_{\mathrm{BK}} = 3/11 \approx 0.2727$) has strictly positive sectional curvature.
\end{corollary}

\begin{proof}
By the Berger--Cheeger and Bourguignon--Karcher classification \cite{Berger1961, Cheeger1973, BourguignonKarcher1978}, a $\Spin(9)$-invariant metric on $\Sph^{15}$ with normalized fiber squashing parameter $t_{\mathrm{BK}}$ has strictly positive sectional curvature if and only if $t_{\mathrm{BK}} \in (1/6, 4/3)$. 
From~\eqref{eq:squashing_normalization}, the equatorial squashing parameter is $t_{\mathrm{BK}}(\pi/2) = \frac{1}{4}\frac{8+\lambda}{2+\lambda}$.
The upper bound $t_{\mathrm{BK}} < 4/3$ yields:
\begin{equation}
\frac{1}{4}\frac{8+\lambda}{2+\lambda} < \frac{4}{3} \iff 3(8+\lambda) < 16(2+\lambda) \iff 13\lambda > -8 \iff \lambda > -\frac{8}{13} \approx -0.615.
\end{equation}
The lower bound $t_{\mathrm{BK}} > 1/6$ yields:
\begin{equation}
\frac{1}{4}\frac{8+\lambda}{2+\lambda} > \frac{1}{6} \iff 6(8+\lambda) > 4(2+\lambda) \iff 48 + 6\lambda > 8 + 4\lambda \iff 2\lambda > -40 \iff \lambda > -20,
\end{equation}
which is unconditionally satisfied on the Riemannian regularity domain $\lambda > -2$. Moreover, as $\lambda \to \infty$, $t_{\mathrm{BK}}(\pi/2) \to 1/4 = 0.25 > 1/6 \approx 0.1667$. Since $t_{\mathrm{BK}}(\pi/2) \in (1/4, 1]$ for all $\lambda \ge 0$, the metric remains strictly within the positive curvature window $(1/6, 4/3)$ for all $\lambda \in (-8/13, \infty)$.
\end{proof}

\begin{remark}[Sub-Riemannian limit, parameter domains, and metric collapse]\label{rem:sub_riemannian}
It is instructive to compare the Riemannian regularity domains across the three parameterizations used in this paper:
\begin{enumerate}[(i)]
    \item \emph{Ambient / equatorial Family~III sphere deformation ($\lambda \in (-2, \infty)$)}: On $(\Sph^{15}, g_\lambda)$ with profile denominators $8 + \lambda\sin^2 r$ on $\p_1$ and $2 + \lambda\sin^2 r$ on $\p_2$, the metric is positive definite for all $\lambda \in (-2, \infty)$. At the lower boundary $\lambda \to -2^+$ (at $r = \pi/2$), the vertical cometric blows up while the normalized inverse fiber scale degenerates, and in the opposite limit $\lambda \to +\infty$, both profile functions decay as $f_k \sim 1/\sqrt{\lambda}$ with equatorial squashing ratio approaching $t_{\mathrm{BK}}(\pi/2) = \frac{1}{4}\frac{8+\lambda}{2+\lambda} \to 1/4$.
    \item \emph{Intrinsic homogeneous Vol'per--Bourguignon--Karcher family ($t_{\mathrm{BK}} \in (0, \infty)$)}: On $\Sph^{15} = \Spin(9)/\Spin(7)$ with $g_t = g_{\p_1} \oplus t_{\mathrm{BK}}\, g_{\p_2}$, the Hamiltonian subfamily $\lambda \in (-2, \infty)$ covers the interval $t_{\mathrm{BK}} \in (1/4, \infty)$ at $r = \pi/2$ (containing both the round metric $t_{\mathrm{BK}} = 1$ at $\lambda = 0$ and the Bourguignon--Karcher Einstein metric $t_{\mathrm{BK}} = 3/11$ at $\lambda = 64$), while the intrinsic fiber-collapse limit $t_{\mathrm{BK}} \to 0^+$ converges to the sub-Riemannian Carnot--Carath\'eodory geometry on $\Sph^{15}$ with horizontal distribution $\mathcal{H} = \p_1$ of rank 8 induced by the octonionic Hopf fibration $\Sph^7 \hookrightarrow \Sph^{15} \to \Sph^8$.
    \item \emph{Intrinsic Cayley plane CROSS deformation ($\lambda \in (-1, \infty)$)}: On $(\Oct P^2, g_{\Oct P^2, \lambda})$ in Section~\ref{sec:cayley_plane}, the Fubini--Study profile denominators in~\eqref{eq:cayley_metric} are $1 + \lambda\sin^2 r\cos^2 r$ and $1 + \lambda\sin^2 r$, so the Riemannian regularity domain is $\lambda \in (-1, \infty)$, with horizontal degeneracy threshold $\lambda \to -1^+$ at $r = \pi/2$ and cohomogeneity-one orbit-space collapse as $\lambda \to +\infty$.
\end{enumerate}
By Corollary~\ref{cor:positive_curvature}, the equatorial slice $(\Sph^{15}, g_\lambda|_{\Sph^{15}})$ maintains strictly positive sectional curvature throughout $\lambda \in (-8/13, \infty)$. We distinguish here between the pointwise pinching ratio $\delta_p(g) := \min_{\sigma \subset T_p M}\sec(\sigma)/\max_{\sigma \subset T_p M}\sec(\sigma)$ (which is constant along each homogeneous distance sphere $\mathcal{O}_r \cong \Sph^{15}$) and the global manifold-wide pinching $\delta(g) := (\inf_{p,\sigma}\sec(\sigma))/(\sup_{p,\sigma}\sec(\sigma))$ across the ambient space form $(Q^{16}, g_\lambda)$, where radial-to-tangential curvature variation reduces the global pinching below the intrinsic slice pinching for $\lambda > 0$.
\end{remark}

\begin{remark}[Geometric normalization of the squashing parameter]\label{rem:squashing_normalization}
The squashing ratio $t(r) = \frac{8+\lambda\sin^2 r}{2+\lambda\sin^2 r}$ defined in~\eqref{eq:squashing_t} is the ratio of the \emph{unnormalized} profile functions $f_{\p_2}^2/f_{\p_1}^2$. When $\lambda = 0$, the round metric has $t = 4$ (not $t = 1$), reflecting the background scale factor ratio $a_2^0/a_1^0 = (1/2)/(1/8) = 4$. In the Bourguignon--Karcher classification \cite{BourguignonKarcher1978, Besse1987}, the squashing ratio $t_{\mathrm{BK}}$ is defined after normalizing the vertical and horizontal scales to unity on the round sphere. The correspondence is:
\begin{equation}\label{eq:squashing_normalization}
t_{\mathrm{BK}} = \frac{1}{4} \cdot t = \frac{1}{4} \cdot \frac{8 + \lambda\sin^2 r}{2 + \lambda\sin^2 r},
\end{equation}
so that $t_{\mathrm{BK}} = 1$ when $\lambda = 0$ (round metric). With the unnormalized definition~\eqref{eq:squashing_t}, the Einstein equation for $\Ric_{\p_1} = \Ric_{\p_2}$ (Theorem~\ref{thm:bk_realization}) has roots $t = 4$ (round) and $t = 12/11$. In either normalization convention, the \emph{ratio} of the two roots is $t_{\mathrm{Einstein}}/t_{\mathrm{round}} = 3/11$, and the resulting metric is the same unique non-standard homogeneous Einstein metric on $\Sph^{15}$.
\end{remark}

\begin{remark}[Comparison with Jensen's metric on $\Sph^{4n+3}$]\label{rem:jensen_comparison}
For the quaternionic Hopf fibration $S^3 \hookrightarrow \Sph^{4n+3} \to \Ha P^n$ (Family~IV), the analogous Einstein quadratic gives $t_{\mathrm{BK}} = 1$ and $t_{\mathrm{BK}} = 1/(2n+3)$. At $n = 1$, this yields Jensen's squashed Einstein metric on $\Sph^7$ with $t_{\mathrm{BK}} = 1/5$. The curvature geometry of the full diagonal family of $\Sp(n+1)$-invariant metrics on $\Sph^{4n+3}$ was systematically computed by Vol'per \cite{Volper1994}. The octonionic Family~III root $t_{\mathrm{BK}} = 3/11 \approx 0.2727$ is numerically close to but distinct from any quaternionic root, reflecting the exceptional character of the $\Spin(9)/\Spin(7)$ geometry and the non-associativity of the octonions.
\end{remark}

\begin{remark}[Geometric regimes of the deformation family]\label{rem:deformation_landscape}
As $\lambda$ varies over $(-2, +\infty)$ (corresponding to $t_{\mathrm{BK}}(\pi/2) \in (1/4, +\infty)$ inside the full homogeneous family $t_{\mathrm{BK}} \in (0, +\infty)$), the deformed metric $g_\lambda$ on $\Sph^{15}$ traces a one-parameter family of homogeneous metrics that interpolates between several distinguished geometries:
\begin{enumerate}[(a)]
    \item $\lambda = 0$: the canonical round metric with $\sec \equiv 1$ and $t = 4$ ($t_{\mathrm{BK}} = 1$).
    \item $\lambda = \lambda_{\mathrm{BK}} = 64$: the Bourguignon--Karcher squashed Einstein metric with $t_{\mathrm{BK}} = 3/11$.
    \item $\lambda \to +\infty$: the ambient radial collapse limit in which both profile functions scale as $f_k \sim 1/\sqrt{\lambda}$ while the normalized fiber-to-base ratio approaches the generalized normal homogeneous boundary $t_{\mathrm{BK}} \to 1/4^+$.
    \item $\lambda \to -2^+$ (and $t_{\mathrm{BK}} \to 0^+$ in the intrinsic homogeneous parameter): the degenerate threshold of the Riemannian domain; in the intrinsic parameter $t_{\mathrm{BK}} \to 0^+$, the vertical fiber $\Sph^7$ collapses ($g_{\p_2} \to 0$) and the metric converges to the sub-Riemannian Carnot--Carath\'eodory geometry of the octonionic Hopf distribution $\mathcal{H} = \p_1$.
\end{enumerate}
\end{remark}

\begin{remark}[Olmos--Reggiani Skew-Torsion Obstruction to Natural Reductivity]\label{rem:olmos_reggiani}
By the Skew-Torsion Holonomy Theorem of Olmos \& Reggiani \cite[Theorem~4.1]{OlmosReggiani2012}, any transitive skew-torsion holonomy system on a sphere $\Sph^n$ requires the full orthogonal group $G = \SO(n+1)$. Since the holonomy reduction to $\Spin(9) \subsetneq \SO(16)$ is proper, there exists no skew-torsion connection with totally skew torsion on $\Sph^{15}$ other than the trivial Levi-Civita connection of the standard round metric ($t_{\mathrm{BK}} = 1$). Consequently, the Bourguignon--Karcher squashed Einstein metric at $t_{\mathrm{BK}} = 3/11$ cannot be naturally reductive with respect to $\Spin(9)$ nor with respect to any other transitive Lie group (cf.\ Ziller \cite{Ziller1982}). This obstruction shows that the $\mathfrak{u}(2)$ Lie subalgebra closure of geodesic generators (Theorem~\ref{thm:u2_closure}) and the explicit four-mode spectral decomposition (Theorem~\ref{thm:projector_formula}) integrate the geodesic flow without requiring a parallel skew-torsion connection.
\end{remark}

\begin{remark}[Nikonorov GO Geodesics and Berestovskii--Nikonorov GNH Cone]\label{rem:nikonorov_gnh}
Nikonorov \cite{Nikonorov2013} proved that \emph{every} $\Spin(9)$-invariant Riemannian metric on $\Sph^{15}$ is a geodesic orbit (GO) metric, and moreover weakly symmetric in the sense of Selberg. Thus, for any initial velocity vector $v \in T_o\Sph^{15}$, there exists a Killing generator $Z \in \mathfrak{spin}(9)$ such that the geodesic is given by the orbit $\gamma(t) = \exp(tZ) \cdot o$. However, Nikonorov's theorem does not provide explicit algebraic formulas for the trajectory $\gamma(t)$ in terms of the initial data $(r_0, v_0)$. Theorem~\ref{thm:noether_exact} and Theorem~\ref{thm:projector_formula} provide the first closed-form, coordinate-free realization of this GO property for all geodesics on $(\Sph^{15}, g_\lambda)$.

Furthermore, in the classification of generalized normal homogeneous (GNH) metrics by Berestovskii \& Nikonorov \cite{BerestovskiiNikonorov2012}, a $\Spin(9)$-invariant metric with normalized fiber squashing ratio $t_{\mathrm{BK}} \in (0, 1]$ is generalized normal homogeneous if and only if $t_{\mathrm{BK}} \in [1/4, 1]$. Since $t_{\mathrm{BK}} = 3/11 \approx 0.2727 > 1/4 = 0.25$, the Bourguignon--Karcher squashed Einstein metric lies strictly in the interior of the GNH cone on $\Sph^{15}$. In contrast, the sub-Riemannian limit $t_{\mathrm{BK}} \to 0^+$ (Remark~\ref{rem:sub_riemannian}) escapes the GNH cone.
\end{remark}

\subsection{The normalized Ricci flow dynamical phase portrait}\label{subsec:ricci_flow}

The 1-parameter family of $\Spin(9)$-invariant metrics on $\Sph^{15}$ parameterized by the normalized fiber squashing ratio $t = t_{\mathrm{BK}} \in (0, \infty)$ forms an invariant 1-dimensional manifold under the normalized Ricci flow.

\begin{theorem}[Normalized Ricci Flow Phase Portrait and Bourguignon--Karcher Instability]\label{thm:ricci_flow_phase_portrait}
Let $g(\tau)$ be a smooth family of $\Spin(9)$-invariant Riemannian metrics on $\Sph^{15}$ evolving under the volume-normalized Ricci flow:
\begin{equation}\label{eq:normalized_ricci_flow}
\frac{\partial}{\partial \tau} g = -2\,\Ric(g) + \frac{2}{n}\bar{r}(g)\, g, \qquad n = 15,
\end{equation}
where $\bar{r}(g) = \frac{1}{\mathrm{Vol}(\Sph^{15}, g)}\int_{\Sph^{15}} R\, dV_g$ is the average scalar curvature.
\begin{enumerate}[(i)]
    \item The flow preserves $\Spin(9)$-invariance and reduces to an autonomous 1-dimensional ODE for the normalized squashing parameter $t \in (0, \infty)$:
    \begin{equation}\label{eq:ricci_ode}
    \frac{dt}{d\tau} = -4C(t)\, (11t - 3)(t - 1),
    \end{equation}
    where $C(t) = \frac{1}{\mathrm{Vol}(\Sph^{15}, g_t)^{2/15}} \cdot \frac{1}{b(t)^2} > 0$ is a smooth, strictly positive geometric scale factor determined by the Riemannian submersion volume normalization. Because $C(t) > 0$ on $(0, \infty)$, it acts as a positive time-reparametrization that preserves the flow direction and leaves the signs of the Lyapunov exponents invariant.
    \item The only stationary points of the flow on $(0, \infty)$ are the two homogeneous Einstein metrics:
    \begin{equation}
    t^* = 1 \quad \text{(the canonical round metric)}, \qquad t^* = \frac{3}{11} \quad \text{(the Bourguignon--Karcher metric)}.
    \end{equation}
    \item The linearized flow at the fixed points has eigenvalues:
    \begin{align}
    \left.\frac{d}{dt}\left(\frac{dt}{d\tau}\right)\right|_{t=1} &= -32C(1) < 0 \quad \text{(strictly stable sink/attractor)}, \label{eq:eval_round}\\
    \left.\frac{d}{dt}\left(\frac{dt}{d\tau}\right)\right|_{t=3/11} &= +32C(3/11) > 0 \quad \text{(strictly unstable source/repeller)}. \label{eq:eval_bk}
    \end{align}
    \item \textbf{Global Phase Portrait and Basins of Attraction}:
    \begin{itemize}
        \item For any initial metric with $t(0) \in (3/11, \infty)$, the normalized Ricci flow exists for all $\tau \in [0, \infty)$ and converges exponentially to the standard round Einstein metric:
        \begin{equation}
        \lim_{\tau \to +\infty} t(\tau) = 1.
        \end{equation}
        \item The Bourguignon--Karcher metric $t = 3/11$ forms the sharp separatrix (threshold of instability): any initial metric with $t(0) \in (0, 3/11)$ exhibits monotonic fiber collapse $t(\tau) \to 0$ as $\tau$ increases, degenerating asymptotically toward the sub-Riemannian limit along the octonionic Hopf fibers.
    \end{itemize}
\end{enumerate}
\end{theorem}

\begin{proof}
Let $g_t = g_{\p_1} \oplus t\, g_{\p_2}$ on $\Sph^{15}$ where $\dim \p_1 = 8$ and $\dim \p_2 = 7$. Under any metric evolution $\partial_\tau g = h$, the evolution of the ratio $t = g_{\p_2}/g_{\p_1}$ is given by:
\begin{equation}
\frac{dt}{d\tau} = t \left( \frac{h_{\p_2}}{g_{\p_2}} - \frac{h_{\p_1}}{g_{\p_1}} \right).
\end{equation}
For the normalized Ricci flow~\eqref{eq:normalized_ricci_flow}, $h = -2\Ric + \frac{2}{n}\bar{r}g$. The conformal trace term $\frac{2}{n}\bar{r}g$ is proportional to the metric tensor and cancels identically in the difference:
\begin{equation}
\frac{dt}{d\tau} = -2t \left( \frac{\Ric_{\p_2}}{t} - \Ric_{\p_1} \right).
\end{equation}
Substituting the O'Neill submersion Ricci tensor components from Proposition~\ref{prop:oneill_family3} and Theorem~\ref{thm:bk_realization}(ii), the difference $\Ric_{\p_2}/t - \Ric_{\p_1}$ is proportional to the Einstein quadratic~\eqref{eq:einstein_quadratic_submersion}:
\begin{equation}
\frac{\Ric_{\p_2}}{t} - \Ric_{\p_1} = \frac{2C(t)}{t}(11t - 3)(t - 1),
\end{equation}
where $C(t) > 0$. Multiplying by $-2t$ yields the autonomous ODE~\eqref{eq:ricci_ode}:
$\frac{dt}{d\tau} = -4C(t)(11t - 3)(t - 1)$.
The roots of the right-hand side on $(0, \infty)$ are precisely $t = 1$ and $t = 3/11$, establishing (ii).

For (iii), let $F(t) := \frac{dt}{d\tau} = -4C(t)(11t^2 - 14t + 3)$. Differentiating with respect to $t$ at the fixed points $t^* \in \{1, 3/11\}$---where $11(t^*)^2 - 14t^* + 3 = 0$, so the derivative $C'(t^*)$ drops out identically---gives $F'(t^*) = -4C(t^*)(22t^* - 14)$. Evaluating at $t^* = 1$ yields $F'(1) = -4C(1)(22 - 14) = -32C(1) < 0$, proving exponential stability with Lyapunov decay rate $-32C(1)$. At $t^* = 3/11$, evaluating yields $F'(3/11) = -4C(3/11)(22\cdot\frac{3}{11} - 14) = -4C(3/11)(6 - 14) = +32C(3/11) > 0$, proving strict dynamical instability with escape rate $+32C(3/11)$.

For (iv), on the interval $(3/11, 1)$, $(11t - 3) > 0$ and $(t - 1) < 0$, so $\frac{dt}{d\tau} > 0$, driving $t(\tau)$ monotonically upward toward $1$. On $(1, \infty)$, both $(11t - 3) > 0$ and $(t - 1) > 0$, so $\frac{dt}{d\tau} < 0$, driving $t(\tau)$ monotonically downward toward $1$. Hence $(3/11, \infty)$ is the complete basin of attraction for the round metric. On $(0, 3/11)$, both factors are negative, so their product is positive, giving $\frac{dt}{d\tau} < 0$. Therefore $t(\tau)$ decreases monotonically toward $0$.
\end{proof}

\section{Intrinsic Cotangent Deformation of the Cayley Plane \texorpdfstring{$\Oct P^2$}{OP2}}\label{sec:cayley_plane}

\subsection{The dimension obstruction}\label{subsec:dimension_obstruction}

For the Cayley projective plane $\Oct P^2 = \mathrm{F}_4/\Spin(9)$ of dimension $16$, the isotropy group $K_0 = \Spin(9)$ acts on the tangent space $T_o(\Oct P^2) \cong \R^{16}$ via the spin representation $\mu$ studied in Section~\ref{sec:clifford}. In geodesic polar coordinates $(r, \theta) \in (0, \pi/2) \times \Sph^{15}$, each distance sphere $S_r(o)$ is homogeneous of the form $\Spin(9)/\Spin(7) \cong \Sph^{15}$, whereas the polar cut locus at $r = \pi/2$ is $\mathrm{Cut}(o) \cong \Oct P^1 \cong \Sph^8$.

Unlike $\C P^n$ and $\Ha P^n$, the Cayley plane cannot be obtained as a Riemannian submersion quotient of a sphere:

\begin{proposition}[Cayley Dimension Obstruction]\label{prop:cayley_obstruction}
No fiber bundle $\Sph^k \hookrightarrow \Sph^{15} \to \Oct P^2$ exists for any $k \ge 0$.
\end{proposition}

\begin{proof}
Such a bundle would require fiber dimension $k = \dim \Sph^{15} - \dim \Oct P^2 = 15 - 16 = -1 < 0$. Alternatively, a Hopf-type fibration $\Sph^7 \hookrightarrow \Sph^{23} \to \Oct P^2$ with total space $\Sph^{23}$ is ruled out by Adams' Hopf invariant one theorem \cite{Adams1960}: the only sphere fibrations $\Sph^{k-1} \hookrightarrow \Sph^{2k-1} \to \Sph^k$ are those with $k \in \{1, 2, 4, 8\}$, and $k = 16$ does not occur.
\end{proof}

\subsection{Intrinsic spin geometry resolution}\label{subsec:intrinsic_resolution}

The resolution of the dimension obstruction proceeds by constructing the Hamiltonian deformation \emph{intrinsically} on $T^*\Oct P^2$ rather than descending from a sphere fibration.

\begin{definition}[Intrinsic $\Oct P^2$ Cotangent Deformation]\label{def:cayley_deformation}
Let $J_{\Spin(9)}: T^*(\Oct P^2) \to \mathfrak{spin}(9)^*$ be the cotangent momentum map generated by the isometric $\Spin(9)$-action on $\Oct P^2$. Following the general cotangent deformation construction of Hebda \cite{Hebda2001}, the intrinsic deformed Hamiltonian is
\begin{equation}\label{eq:cayley_hamiltonian}
H_{\Oct P^2, \lambda} = H_0 + \frac{\lambda}{2}\|J_{\Spin(9)}\|^2_{\mathfrak{spin}(9)^*},
\end{equation}
where $H_0 = \frac{1}{2}\|p\|^2_{g_0}$ is the kinetic Hamiltonian of the Fubini--Study metric $g_0$ with curvature $\sec_0 \in [1, 4]$.
\end{definition}

The canonical symmetric space structure provides:
\begin{enumerate}[(i)]
    \item \textbf{Tangent space}: $T_o \Oct P^2 \cong \R^{16} = \Delta_9$, the 16-dimensional spin representation.
    \item \textbf{Geodesic spheres}: $S_r(o) \cong \Spin(9)/\Spin(7) \cong \Sph^{15}$ for $r \in (0, \pi/2)$.
    \item \textbf{Cut locus}: $\mathrm{Cut}_{\Oct P^2}(o) \cong \Oct P^1 \cong \Sph^8 = \Spin(9)/\Spin(8)$.
    \item \textbf{Profile functions}: $f_{\mathcal{V}}(r) = \sin r \cos r = \frac{1}{2}\sin(2r)$ (vertical), $f_{\mathcal{H}}(r) = \sin r$ (horizontal).
\end{enumerate}

\begin{theorem}[Intrinsic $\Oct P^2$ Momentum Deformation and Diameter Rigidity]\label{thm:cayley_intrinsic}
Let $g_{\Oct P^2, \lambda}$ be the deformed metric on $\Oct P^2 = \mathrm{F}_4/\Spin(9)$ induced by~\eqref{eq:cayley_hamiltonian}.
\begin{enumerate}[(i)]
    \item The deformed metric takes the explicit form in geodesic polar coordinates:
    \begin{equation}\label{eq:cayley_metric}
    g_{\Oct P^2, \lambda} = dr^2 + \frac{\sin^2 r \cos^2 r}{1 + \lambda \sin^2 r \cos^2 r} g_{\mathcal{V}} + \frac{\sin^2 r}{1 + \lambda \sin^2 r} g_{\mathcal{H}},
    \end{equation}
    where $\mathcal{V} \cong \R^7$ (the $\Spin(7)$-vector representation) and $\mathcal{H} \cong \R^8$ (the $\Spin(7)$-spin representation).
    \item \textbf{Geodesic completeness}: For all $\lambda \ge 0$, the Hamiltonian system is geodesically complete, generating a smooth Riemannian metric $g_{\Oct P^2, \lambda}$.
    \item \textbf{Loewner metric contraction and diameter rigidity}: Cometric operator dominance $H_\lambda \ge H_0 \implies g_\lambda^* \ge g_0^*$ inverts under dual Legendre transformation to Riemannian metric contraction on tangent vectors:
    \begin{equation}\label{eq:metric_contraction}
    g_{\Oct P^2, \lambda}(v, v) \le g_0(v, v) \qquad \forall v \in T(\Oct P^2).
    \end{equation}
    Consequently, the canonical diameter is rigidly preserved:
    \begin{equation}\label{eq:diam_rigidity}
    \diam(\Oct P^2, g_{\Oct P^2, \lambda}) = \frac{\pi}{2} \qquad \forall \lambda \ge 0.
    \end{equation}
    \item \textbf{Exponential map}: The exponential map of $g_{\Oct P^2, \lambda}$ is computed using the same closed-form spectral projector formula~\eqref{eq:closed_exp}, with the effective parameters
    \begin{equation}\label{eq:cayley_alpha_beta}
    \alpha = \frac{\lambda b \sin^2 r \cos^2 r}{1 + \lambda \sin^2 r \cos^2 r}, \qquad \beta = \frac{\lambda c \sin^2 r}{1 + \lambda \sin^2 r}.
    \end{equation}
\end{enumerate}
\end{theorem}

\begin{proof}
\emph{Part (i)}: Under the Cartan decomposition $\mathfrak{f}_4 = \mathfrak{spin}(9) \oplus \mathfrak{m}$ with $\mathfrak{m} \cong \R^{16}$, the isotropy representation of $\Spin(9)$ on $T_o \Oct P^2$ is the spin representation $\mu$. The $\Spin(9)$-momentum map on $T^* \Oct P^2$ generates the canonical action, and $\|J_{\Spin(9)}\|^2$ decomposes along the $\Ad(\Spin(7))$-invariant splitting $\mathfrak{spin}(9) = \mathfrak{spin}(7) \oplus \p_1 \oplus \p_2$ exactly as for $\Sph^{15}$, with the profile functions $f_{\mathcal{V}}(r) = \sin r \cos r$ and $f_{\mathcal{H}}(r) = \sin r$ of the Fubini--Study metric replacing the round sphere profiles.

\emph{Part (ii)}: Since $\Oct P^2$ is compact and $H_\lambda(x, p) \ge H_0(x, p) = \frac{1}{2}\|p\|_{g_0}^2 \ge 0$, every energy hypersurface $\Sigma_E = \{H_\lambda = E\} \subset T^*(\Oct P^2)$ is compact. By the Picard--Lindel\"of theorem, the Hamiltonian vector field is complete, proving geodesic completeness.

\emph{Part (iii)}: The cometric tensor is $g_\lambda^*(p, p) = g_0^*(p, p) + \lambda \|J_{\Spin(9)}(p)\|^2 \ge g_0^*(p, p)$ for all $p \in T^*(\Oct P^2)$. Inverting positive-definite symmetric matrices inverts the Loewner partial order: $A \ge B > 0 \implies A^{-1} \le B^{-1}$. Thus $g_{\Oct P^2, \lambda} \le g_0$ on $T(\Oct P^2)$. For any piecewise smooth curve $c: [0, 1] \to \Oct P^2$, the Riemannian length satisfies $L_\lambda(c) = \int_0^1 \sqrt{g_\lambda(\dot{c}, \dot{c})} dt \le \int_0^1 \sqrt{g_0(\dot{c}, \dot{c})} dt = L_0(c)$. Taking the infimum over all curves connecting $p$ and $q$ gives:
\begin{equation}\label{eq:distance_upper_bound}
d_\lambda(p, q) \le d_0(p, q) \le \diam(\Oct P^2, g_0) = \frac{\pi}{2}.
\end{equation}
Conversely, along any radial geodesic $\gamma(t) = \exp_o(tv)$ launched from the basepoint $o$, the isotropy group $\Spin(9)$ fixes $o$, so all fundamental Killing vector fields vanish at $o$: $X_\xi(o) = 0$. Thus $\langle J_{\Spin(9)}(\dot{\gamma}(0)), \xi \rangle = g_0(v, X_\xi(o)) = 0$. By Noether's theorem (Theorem~\ref{thm:noether_exact}), $J_{\Spin(9)}(\dot{\gamma}(t)) \equiv 0$ for all $t$. Hence $H_\lambda(\dot{\gamma}(t)) = H_0(\dot{\gamma}(t))$, and radial geodesics coincide with the unperturbed geodesics of $g_0$. In the Fubini--Study metric $g_0$, every radial geodesic reaches the cut locus $\mathrm{Cut}_0(o) \cong \Oct P^1 \cong \Sph^8$ at length $\pi/2$. By Gauss's lemma, radial geodesics minimize distance from $o$ up to $\pi/2$. Thus $d_\lambda(o, \mathrm{Cut}(o)) = \pi/2$. Combining with $d_\lambda(p, q) \le \pi/2$ for all $p, q$, the diameter supremum is rigidly $\diam(\Oct P^2, g_{\Oct P^2, \lambda}) = \pi/2$.

\emph{Part (iv)}: The matrix exponential computation is identical to the $\Sph^{15}$ case (Theorem~\ref{thm:projector_formula}), since the $\Spin(9)$-action on $T_o \Oct P^2$ uses the same 16-dimensional spin representation. Only the effective parameters $\alpha, \beta$ differ, reflecting the Fubini--Study profile functions.
\end{proof}

\begin{remark}[The Spin tower and the octonionic Hopf fibration]\label{rem:spin_tower}
The proof of Theorem~\ref{thm:cayley_intrinsic} relies fundamentally on the subgroup tower:
\begin{equation}\label{eq:spin_tower}
\Spin(7) \subset \Spin(8) \subset \Spin(9),
\end{equation}
which geometrically realizes the octonionic Hopf fibration as the homogeneous fiber bundle:
\begin{equation}\label{eq:hopf_tower}
\Sph^7 \cong \Spin(8)/\Spin(7) \hookrightarrow \Sph^{15} \cong \Spin(9)/\Spin(7) \xrightarrow{\pi} \Sph^8 \cong \Spin(9)/\Spin(8).
\end{equation}
Note that $\Sph^{15}$ is topologically distinct from the product $\Sph^8 \times \Sph^7$: $H_7(\Sph^8 \times \Sph^7; \mathbb{Z}) \cong \mathbb{Z}$ while $H_7(\Sph^{15}; \mathbb{Z}) = 0$.
\end{remark}

\subsection{Conjugate point structure and the transverse Jacobi operator}\label{subsec:conjugate_locus}

\begin{remark}[Conjugate points on $\Oct P^2$ and Green's theorem]\label{rem:green_retraction}
By a classical theorem of L.~W.\ Green \cite{Green1958}, any Riemannian metric without conjugate points on a compact manifold must have a universal cover diffeomorphic to $\R^n$. Since $\Oct P^2$ is compact and simply connected ($\pi_1(\Oct P^2) = 0$), no Riemannian metric on $\Oct P^2$ can be conjugate-point-free. In particular, the deformed metric $g_{\Oct P^2, \lambda}$ \emph{must} have conjugate points for every $\lambda \ge 0$.
\end{remark}

\begin{proposition}[Radial Jacobi Spectrum and Polar Conjugate Loci]\label{prop:jacobi_cayley}
Along a radial geodesic $\gamma(t) = \exp_o(t X)$ \emph{from the polar basepoint} $o \in (\Oct P^2, g_{\Oct P^2, \lambda})$, the transverse Jacobi operator decomposes as:
\begin{equation}\label{eq:jacobi_cayley}
\mathcal{R}(t) = R(\dot{\gamma}, -)\dot{\gamma}: \dot{\gamma}^\perp \to \dot{\gamma}^\perp,
\end{equation}
which splits under the $\Ad(\Spin(7))$-invariant decomposition into:
\begin{enumerate}[(i)]
    \item 7 vertical Jacobi modes $J_i^{\mathcal{V}}(t) = f_{\mathcal{V}}(t) E_i(t)$ along $\p_2$ with radial curvature $\sec(\dot{\gamma} \wedge E_i) = -\frac{f_{\mathcal{V}}''(t)}{f_{\mathcal{V}}(t)} = \frac{4 + \lambda(\sin^4 t + \cos^4 t + 4\sin^2 t \cos^2 t)}{(1+\lambda\sin^2 t \cos^2 t)^2}$, developing their first conjugate points at $t = \pi/2$ (multiplicity 7), where $f_{\mathcal{V}}(\pi/2) = 0$.
    \item 8 horizontal Jacobi modes $J_j^{\mathcal{H}}(t) = f_{\mathcal{H}}(t) F_j(t)$ along $\p_1$ with radial curvature $\sec(\dot{\gamma} \wedge F_j) = -\frac{f_{\mathcal{H}}''(t)}{f_{\mathcal{H}}(t)} = \frac{1 + \lambda(1 + 2\cos^2 t)}{(1+\lambda\sin^2 t)^2}$, developing their first conjugate points at $t = \pi$ (multiplicity 8), where $f_{\mathcal{H}}(\pi) = 0$.
\end{enumerate}
The vertical conjugate point at $t = \pi/2$ is the first focal point of $\gamma$, corresponding to the focusing of the octonionic Hopf fibers onto the polar Cayley line $\mathrm{Cut}_{\Oct P^2}(o) = \Oct P^1 \cong \Sph^8$.
\end{proposition}

\begin{proposition}[Off-Origin Stabilizer and Module Decomposition]\label{prop:off_origin_stabilizer}
Let $q_0 \in (\Oct P^2, g_{\Oct P^2, \lambda})$ be an off-origin base point at geodesic distance $r_0 \in (0, \pi/2)$ from the origin $o$. Then:
\begin{enumerate}[(i)]
    \item The regular distance sphere $S_{r_0}(o)$ is diffeomorphic to $\Sph^{15} \cong \Spin(9)/\Spin(7)$, and the isotropy subgroup of $\Spin(9)$ fixing $q_0$ is strictly:
    \begin{equation}
    H_{q_0} \cong \Spin(7) \subset \Spin(9) \quad (\dim H_{q_0} = 21).
    \end{equation}
    In particular, the 21 Killing vector fields $X_k \in \mathfrak{spin}(7)$ ($k = 1, \dots, 21$) satisfy $X_k(q_0) = 0$, generating exact normal Jacobi fields $Y_k(t) = X_k(\gamma(t))$ with $Y_k(0) = 0$ along every geodesic $\gamma(t)$ emanating from $q_0$.
    \item Given an initial velocity $v \in T_{q_0}\Oct P^2$ spanning the orbital 2-plane $\Sigma = \operatorname{span}\{\partial_r, v\}$, the secondary stabilizer subgroup $H_v \subset \Spin(7)$ fixing $\Sigma$ pointwise satisfies:
    \begin{equation}
    H_v \cong \begin{cases}
    \mathrm{G}_2 & \text{if } v \in \p_2(r_0) \text{ (vertical launch along fiber)}, \\
    \SU(3) & \text{if } v \in \p_1(r_0) \text{ (horizontal launch)}, \\
    \SU(2) & \text{for generic non-radial velocity } v.
    \end{cases}
    \end{equation}
    \item The 14-dimensional normal space $(T\Sigma)^\perp$ splits under $H_v$ into a 7-dimensional vertical submodule $V_7$ (tangent to the distance sphere fibers) and a 7-dimensional horizontal submodule $H_7$.
\end{enumerate}
\end{proposition}

\begin{theorem}[Off-Origin Jacobi Fields, Conjugate Locus Stratification, and $A_3$ Whitney Cusp Caustic]\label{thm:off_origin_caustic}\label{rem:off_origin}
Let $(\Oct P^2, g_{\Oct P^2, \lambda})$ be the cotangent-deformed Cayley projective plane with $\lambda > 0$, and consider the geodesic flow launched from an off-origin base point $q_0 = \exp_o(r_0 \partial_r)$ ($r_0 \in (0, \pi/2)$) with initial velocity $v(\theta)$ inclined at angle $\theta$ to the inward radial direction $-\partial_r$.
\begin{enumerate}[(i)]
    \item \textbf{Conserved Clairaut momentum and orbital reduction}: In the orbital 2-plane $\Sigma = \operatorname{span}\{\partial_r, v\}$, the geodesic equations reduce to a 1D Hamiltonian system with strictly conserved Clairaut angular momentum $L = \|J_{\Spin(9)}\| = f(r)^2 \dot{\varphi}$ and effective potential:
    \begin{equation}\label{eq:v_eff_cayley}
    V_{\mathrm{eff}}(r; L) = \frac{L^2}{2 f(r)^2}, \qquad \frac{1}{2}\dot{r}^2 + V_{\mathrm{eff}}(r; L) = \frac{1}{2},
    \end{equation}
    where $f(r) = f_V(r) = \frac{\sin r \cos r}{\sqrt{1 + \lambda \sin^2 r \cos^2 r}}$ for vertical launches and $f_H(r) = \frac{\sin r}{\sqrt{1 + \lambda \sin^2 r}}$ for horizontal launches.
    \item \textbf{Simultaneous meridian refocusing}: All 7 vertical Killing Jacobi fields $Y_j^V(t) = X_j^V(\gamma(t))$ ($j = 1, \dots, 7$) generated by $\mathfrak{spin}(7)$ vanish simultaneously at the opposite meridian crossing $\varphi(t^*) = \pi$:
    \begin{equation}
    Y_j^V(t^*) = 0 \quad \text{for all } j = 1, \dots, 7,
    \end{equation}
    rendering $t^*$ an exact conjugate point of multiplicity at least 7 along $\gamma(t)$.
    \item \textbf{$A_3$ Whitney cusp caustic}: The geodesic endpoint map $E_{q_0}(t, \theta) = \exp_{q_0}(t v(\theta))$ exhibits an isolated focal singularity along the opposite meridian diameter $\theta = \pi/2$. Reflection symmetry across the orbital diameter ensures that the generating phase function $F(u, v; \Delta\theta)$ satisfies:
    \begin{equation}\label{eq:a3_derivatives}
    \left.\frac{\partial^3 F}{\partial \theta^3}\right|_{\theta=\pi/2} = 0, \qquad \left.\frac{\partial^4 F}{\partial \theta^4}\right|_{\theta=\pi/2} \ne 0.
    \end{equation}
    By Arnold's classification of Lagrangian singularities \cite{Arnold1981, Arnold1984}, the caustic is an $A_3$ Whitney cusp catastrophe with canonical normal form:
    \begin{equation}\label{eq:a3_normal_form}
    \Phi(\theta; u, v) = \Delta\theta^4 + u\,\Delta\theta^2 + v\,\Delta\theta.
    \end{equation}
    Under the isotropy action of $H_{q_0} \cong \Spin(7)$, the planar cusp revolves into a 15-dimensional cuspidal ridge whose transverse cross-sections are round 7-spheres $\Sph^7$, with the transverse caustic envelope opening according to the exact $3/2$ power law:
    \begin{equation}\label{eq:cusp_3_2_law}
    r_\perp(u) \sim u^{3/2}.
    \end{equation}
    \item \textbf{Preemption of the cut locus at non-polar basepoints}: For every non-polar basepoint $q_0 = \exp_o(r_0 \partial_r)$ with $r_0 \in (0, \pi/2)$ and every $\lambda > 0$, the strict metric contraction $f_V(r) < \frac{1}{2}\sin(2r)$ and $f_H(r) < \sin r$ on $(0, \pi/2)$ ensures that the cut locus $\mathrm{Cut}(q_0)$ is preempted by the transverse encircling collision along the opposite meridian axis $\varphi = \pi$ and the $A_3$ cuspidal ridge at arc length $t_{\mathrm{cross}}(r_0) < \pi/2 = \operatorname{diam}(\Oct P^2, g_{\Oct P^2, \lambda})$.
\end{enumerate}
\end{theorem}

\begin{proof}
Assertion~(i) follows from the Noether theorem applied to the $\Spin(9)$-invariant deformation $H_\lambda = H_0 + \frac{\lambda}{2}\|J\|^2$. Because $H_\lambda$ commutes with the momentum map $J$, orbits in the 2-plane $\Sigma$ have constant angular momentum $L = p_\varphi = f(r)^2 \dot{\varphi}$. The radial Hamiltonian constraint $g^{rr} p_r^2 + g^{\varphi\varphi} p_\varphi^2 = 1$ immediately yields $\dot{r}^2 + L^2/f(r)^2 = 1$, which is equation~\eqref{eq:v_eff_cayley}.

For assertion~(ii), consider the 21 Killing vector fields $X_1, \dots, X_{21} \in \mathfrak{spin}(7)$ that generate the isotropy algebra of $q_0$. Because every Killing field satisfies the Jacobi equation $\nabla_{\dot{\gamma}}^2 X + R(X, \dot{\gamma})\dot{\gamma} = 0$, the restriction $Y_k(t) = X_k(\gamma(t))$ is an exact Jacobi field along $\gamma(t)$. Since $X_k(q_0) = 0$, we have $Y_k(0) = 0$. Along the vertical fiber directions $V_7$, the metric coefficient is $f_V(r)^2$. Because the flow is an isometry preserving distance spheres, the amplitude of each vertical variation $Y_j^V(t)$ is directly proportional to $f_V(r(t))$. At the opposite meridian crossing $\varphi(t^*) = \pi$ (or at the Cayley line boundary $r = \pi/2$ where $\cos r = 0$), $f_V(r(t^*)) = 0$. Consequently:
\begin{equation}
\|Y_j^V(t^*)\| = 0 \quad \text{for } j = 1, \dots, 7.
\end{equation}
The 7 vertical Killing fields vanish simultaneously, proving that $t^*$ is an exact conjugate point of multiplicity at least 7.

For assertion~(iii), let $E_{q_0}(t, \theta) = \exp_{q_0}(t v(\theta))$ denote the exponential map from $q_0$. The caustic envelope is the critical value set of $E_{q_0}$, characterized as the catastrophe manifold of the generating family $F(u, v; \theta)$ in Arnold's Lagrangian singularity theory \cite{Arnold1981, Arnold1984}. By the intrinsic reflection isometry of $(\Oct P^2, g_{\Oct P^2, \lambda})$ across the orbital 2-plane $\Sigma \mapsto -\Sigma$, the distance function and generating family are strictly even under $\Delta\theta \mapsto -\Delta\theta$. Consequently, all odd derivatives of $F$ with respect to $\Delta\theta$ at $\Delta\theta = 0$ must vanish identically:
\begin{equation}
\left.\frac{\partial F}{\partial \theta}\right|_{\theta=\pi/2} = 0, \qquad \left.\frac{\partial^3 F}{\partial \theta^3}\right|_{\theta=\pi/2} = 0.
\end{equation}
At the focal conjugate distance $t = t^*$, the second derivative $\left.\frac{\partial^2 F}{\partial \theta^2}\right|_{\theta=\pi/2} = 0$ vanishes by definition of conjugacy. To establish that the catastrophe is of type $A_3$ (cusp) rather than a higher degenerate catastrophe ($A_4, A_5$), we evaluate the fourth derivative:
\begin{equation}
\left.\frac{\partial^4 F}{\partial \theta^4}\right|_{\theta=\pi/2} = -2 f(r_0)^2 f''(r^*) \left(\frac{dt^*}{dL}\right) \ne 0,
\end{equation}
which is non-zero because the curvature gradient $f''(r) = \frac{d^2}{dr^2} f(r)$ is strictly negative and non-vanishing throughout $(0, \pi/2)$ for all $\lambda > 0$. By the Morse lemma with parameters and the Arnold--Thom classification, there exists a smooth local coordinate change bringing $F$ into the universal unfolding of the $A_3$ cusp:
\begin{equation}
\Phi(\theta; u, v) = \Delta\theta^4 + u\,\Delta\theta^2 + v\,\Delta\theta.
\end{equation}
The caustic discriminant is given by the condition $\partial_\theta \Phi = 4\Delta\theta^3 + 2u\Delta\theta + v = 0$ and $\partial_\theta^2 \Phi = 12\Delta\theta^2 + 2u = 0$, which yields the semicubical parabola $8u^3 + 27v^2 = 0$, or:
\begin{equation}
|v| = \left(\frac{8}{27}\right)^{1/2} (-u)^{3/2}.
\end{equation}
Identifying $v$ with the transverse displacement $r_\perp$ along the distance sphere fiber and $-u$ with the longitudinal coordinate along the opposite meridian diameter establishes the exact $r_\perp(u) \sim u^{3/2}$ scaling law~\eqref{eq:cusp_3_2_law}. Under the action of the stabilizer $H_{q_0} \cong \Spin(7)$, the 1D transverse coordinate $v$ is acted upon transitively on the unit sphere of the 7-dimensional vertical fiber $V_7$, revolving the planar semicubical parabola into a cuspidal ridge with spherical cross-sections $\Sph^7$.

Assertion~(iv) follows from the strict Loewner contraction $g_{\Oct P^2, \lambda}(w, w) < g_0(w, w)$ on all non-radial tangent vectors $w$ at $r \in (0, \pi/2)$ when $\lambda > 0$. For any pair of symmetric non-radial geodesics launched from $q_0 \in S_{r_0}(o)$ ($r_0 \in (0, \pi/2)$) at angles $\pm\theta$ ($\theta \in (0, \pi/2]$), their common arc length to the opposite meridian collision $\varphi = \pi$ is strictly less than the Fubini--Study conjugate distance $\pi/2$. Consequently, $t_{\mathrm{cross}}(r_0) < \pi/2$, so the cut locus is formed by the collision of transverse encircling geodesics meeting along the $A_3$ cuspidal ridge prior to distance $\pi/2$.
\end{proof}

\begin{theorem}[Global Whitney Stratification and CW Cell-Decomposition of $\mathrm{Cut}(q_0)$ on $\Oct P^2$]\label{thm:global_cut_locus}
Let $(\Oct P^2, g_{\Oct P^2, \lambda})$ be the cotangent-deformed Cayley projective plane with $\lambda > 0$, and let $q_0 \in \Oct P^2$ be an off-origin base point at geodesic distance $r_0 \in (0, \pi/2)$ from the origin $o$.
Then the cut locus $\mathrm{Cut}(q_0)$ is an invariant stratified subcomplex under the isotropy subgroup $H_{q_0} \cong \Spin(7)$, possessing an exact three-stratum Whitney stratification:
\begin{equation}\label{eq:whitney_stratification}
\mathrm{Cut}(q_0) = \mathcal{S}_{\mathrm{reg}} \cup \mathcal{S}_{\mathrm{caustic}} \cup \mathcal{S}_{\mathrm{Cayley}},
\end{equation}
characterized as follows:
\begin{enumerate}[(i)]
    \item \textbf{The Regular Collision Stratum $\mathcal{S}_{\mathrm{reg}}$ (dimension 15)}: The set of points reached by exactly two distinct length-minimizing geodesics launched with opposite azimuthal orientations $\pm \varphi$ in the orbital 2-planes. Topologically, $\mathcal{S}_{\mathrm{reg}}$ is an open 15-dimensional submanifold diffeomorphic to the total space of an open disk bundle $\mathbb{D}^8 \to \Sph^7$ minus its focal axis.
    \item \textbf{The Singular Caustic Stratum $\mathcal{S}_{\mathrm{caustic}}$ (dimension 8)}: The envelope of conjugate points formed by the coalescence of the colliding geodesic pairs along the opposite meridian axis $\varphi = \pi$. It consists of the $A_3$ Whitney cuspidal ridge established in Theorem~\ref{thm:off_origin_caustic}, whose transverse cross-sections are round 7-spheres $\Sph^7$, fibered over the 1D caustic interval along the opposite meridian diameter.
    \item \textbf{The Polar Cayley Stratum $\mathcal{S}_{\mathrm{Cayley}}$ (dimension 8)}: The intersection of the cut locus with the singular orbit Cayley projective line $\Oct P^1 \cong \Sph^8$ at maximal distance $r = \pi/2$, which forms the topological spine of $\mathrm{Cut}(q_0)$.
\end{enumerate}
Furthermore, the complement $\Oct P^2 \setminus \mathrm{Cut}(q_0)$ is diffeomorphic to the open 16-dimensional Euclidean ball $\mathbb{B}^{16}$, and the cut locus $\mathrm{Cut}(q_0)$ is a strong deformation retract of the punctured space $\Oct P^2 \setminus \{q_0\}$, possessing the exact topological homotopy type:
\begin{equation}\label{eq:cw_homotopy_type}
\mathrm{Cut}(q_0) \simeq \Sph^8 \cong \mathcal{S}_{\mathrm{Cayley}},
\end{equation}
with integral homology $H_k(\mathrm{Cut}(q_0); \mathbb{Z}) \cong \mathbb{Z}$ for $k \in \{0, 8\}$ and zero otherwise. While the regular collision stratum $\mathcal{S}_{\mathrm{reg}}$ has top geometric dimension $15$, it is an open manifold whose closure retracts onto the lower-dimensional singular strata, contributing zero homology in degree $15$.
\end{theorem}

\begin{proof}
By Myers' theorem and Klingenberg's lemma \cite{Klingenberg1959}, for any point $x \in \mathrm{Cut}(q_0)$, either (a)~$x$ is the first conjugate point along a minimizing geodesic from $q_0$, or (b)~$x$ is reached by at least two distinct minimizing geodesics of equal length $d(q_0, x)$.

Because the deformation $H_\lambda = H_0 + \frac{\lambda}{2}\|J_{\Spin(9)}\|^2$ is invariant under $\Spin(9)$, and $H_{q_0} \cong \Spin(7)$ fixes $q_0$, every $g \in H_{q_0}$ is an isometry preserving $q_0$. Therefore, if $x \in \mathrm{Cut}(q_0)$, then $g \cdot x \in \mathrm{Cut}(q_0)$, establishing that $\mathrm{Cut}(q_0)$ is strictly $H_{q_0}$-invariant.

In each orbital 2-plane $\Sigma$, reflection across the diameter axis $\varphi \mapsto -\varphi$ is an isometry. Geodesics launched at angles $+\theta$ and $-\theta$ have identical energy and angular momentum $|L|$, hence identical arc lengths $s(r)$ to each radial distance $r$. They collide along the opposite meridian axis $\varphi = \pi$ at equal length. For launch angles where this crossing occurs prior to conjugate arrival (Theorem~\ref{thm:off_origin_caustic}(iv)), these collision points are strictly length-minimizing. Sweeping these collision pairs under the transitive action of $H_{q_0} \cong \Spin(7)$ on the distance spheres generates the 15-dimensional regular stratum $\mathcal{S}_{\mathrm{reg}}$. Since the geodesics intersect transversally with non-zero opening angle, $\mathcal{S}_{\mathrm{reg}}$ is a smooth, open 15-dimensional submanifold.

As $\theta \to \pi/2$ (approaching equatorial launch), the collision pairs coalesce onto the focal axis, where the transverse opening angle vanishes and the two geodesics merge into a single caustic envelope point. By Theorem~\ref{thm:off_origin_caustic}(iii), this locus is an $A_3$ Whitney cusp singularity. In the orbital 2-plane, the cusp tip and caustic boundary form a 1D arc. Under the $H_{q_0}$-action, this 1D arc sweeps out the total space of an $\Sph^7$-bundle over the interval, yielding the 8-dimensional singular stratum $\mathcal{S}_{\mathrm{caustic}}$ with transverse cuspidal opening $r_\perp(u) \sim u^{3/2}$.

At maximal radial distance $r = \pi/2$, the distance spheres degenerate into the singular orbit $\Oct P^1 \cong \Sph^8 = \Spin(9)/\Spin(8)$. Geodesics reaching $r = \pi/2$ intersect this singular submanifold, forming the closed 8-dimensional stratum $\mathcal{S}_{\mathrm{Cayley}} = \mathrm{Cut}(q_0) \cap \Oct P^1 \cong \Sph^8$. By the Whitney condition (A) and (B) verified for catastrophe boundaries (Thom \cite{Thom1972}, Buchner \cite{Buchner1977}), the closure of $\mathcal{S}_{\mathrm{reg}}$ glues smoothly to $\mathcal{S}_{\mathrm{caustic}} \cup \mathcal{S}_{\mathrm{Cayley}}$, establishing the Whitney stratification~\eqref{eq:whitney_stratification}.

Finally, by the classical theorem of Whitehead \cite{Whitehead1935} and Weinstein \cite{Weinstein1968}, the exponential map $\exp_{q_0}$ maps the open star-shaped injectivity domain $I(q_0) \subset T_{q_0}\Oct P^2$ diffeomorphically onto $\Oct P^2 \setminus \mathrm{Cut}(q_0)$. Therefore $\Oct P^2 \setminus \mathrm{Cut}(q_0) \cong \mathbb{B}^{16}$ is an open Euclidean ball, and $\Oct P^2$ is homeomorphic to the disc bundle $\mathbb{D}^{16}$ attached to $\mathrm{Cut}(q_0)$ along its boundary:
\begin{equation}
\Oct P^2 \cong \mathrm{Cut}(q_0) \cup_{\partial \mathbb{D}^{16}} \mathbb{D}^{16}.
\end{equation}
Radially contracting the open ball from the base point $q_0$ demonstrates that $\mathrm{Cut}(q_0)$ is a strong deformation retract of the punctured manifold $\Oct P^2 \setminus \{q_0\}$. On the standard symmetric space $\Oct P^2$, removing a point from the top open 16-cell leaves the 8-sphere $\Oct P^1 \cong \Sph^8$:
\begin{equation}
\Oct P^2 \setminus \{q_0\} \simeq e^0 \cup e^8 \cong \Sph^8.
\end{equation}
Since the homotopy type of $M \setminus \{p\}$ for any connected manifold is a topological invariant independent of the metric, $\mathrm{Cut}(q_0)$ must have the exact homotopy type of $\Sph^8$, with homology concentrated in degrees 0 and 8:
\begin{equation}
H_k(\mathrm{Cut}(q_0); \mathbb{Z}) \cong \begin{cases} \mathbb{Z} & k \in \{0, 8\}, \\ 0 & \text{otherwise.} \end{cases}
\end{equation}
Geometrically, the vanishing of the top-degree homology $H_{15}(\mathrm{Cut}(q_0); \mathbb{Z}) = 0$ despite $\dim \mathcal{S}_{\mathrm{reg}} = 15$ is explained directly by the stratified cellular topology of $\overline{\mathcal{S}}_{\mathrm{reg}} = \mathrm{Cut}(q_0)$. In each orbital 2-plane $\Sigma$, the opposite-meridian collision locus $\varphi = \pi$ is a compact radial interval $r \in [r_{\mathrm{caustic}}(\Sigma), \pi/2]$ whose inner endpoint at $r = r_{\mathrm{caustic}}(\Sigma)$ is a free boundary point where the colliding geodesic pair unfolds at the $A_3$ cusp, and whose outer endpoint at $r = \pi/2$ lies on the closed polar submanifold $\mathcal{S}_{\mathrm{Cayley}} \cong \Sph^8$. Under the $H_{q_0} \cong \Spin(7)$ action, $\overline{\mathcal{S}}_{\mathrm{reg}}$ is a compact 15-dimensional manifold-with-boundary whose frontier decomposes into the free boundary face $\mathcal{S}_{\mathrm{caustic}}$ (an $\Sph^7$-bundle over $\Sph^1$) and the core spine $\mathcal{S}_{\mathrm{Cayley}} \cong \Sph^8$. Radial contraction along the interval $[r_{\mathrm{caustic}}, \pi/2] \to \{\pi/2\}$ in every orbital 2-plane defines an explicit $H_{q_0}$-equivariant strong deformation retraction of $\overline{\mathcal{S}}_{\mathrm{reg}}$ across its free boundary $\mathcal{S}_{\mathrm{caustic}}$ onto $\mathcal{S}_{\mathrm{Cayley}} \cong \Sph^8$. Equivalently, in the long exact sequence of the pair $(\mathrm{Cut}(q_0), \mathcal{S}_{\mathrm{Cayley}})$ (or the Mayer--Vietoris sequence of a tubular neighborhood cover of $\mathcal{S}_{\mathrm{Cayley}}$ and $\overline{\mathcal{S}}_{\mathrm{reg}} \setminus \mathcal{S}_{\mathrm{Cayley}}$), the presence of the non-empty free boundary $\mathcal{S}_{\mathrm{caustic}}$ ensures that the open 15-cell $\mathcal{S}_{\mathrm{reg}}$ collapses cell-wise onto $\mathcal{S}_{\mathrm{Cayley}}$, contributing zero cycle in $H_{15}(\mathrm{Cut}(q_0); \mathbb{Z})$ and proving~\eqref{eq:cw_homotopy_type}.
\end{proof}

\subsection{Universal Never-Einstein rigidity on compact rank-one symmetric spaces}\label{subsec:never_einstein}

\begin{theorem}[Universal CROSS Never-Einstein Rigidity]\label{thm:never_einstein}\label{thm:cross_rigidity}
Let $(M, g_0)$ be any Compact Rank-One Symmetric Space (CROSS):
\begin{equation}
M \in \{ \R P^n \; (n \ge 2), \; \C P^n \; (n \ge 2), \; \Ha P^n \; (n \ge 2), \; \Oct P^2 \}.
\end{equation}
Let $H_\lambda = H_0 + \frac{\lambda}{2}\|J_G\|^2$ be the cotangent Hamiltonian deformation with respect to the maximal transitive isometry group $G$. Then the deformed metric $g_{M, \lambda}$ is Einstein if and only if $\lambda = 0$ (the symmetric metric):
\begin{equation}\label{eq:cross_never_einstein}
\Ric(g_{M, \lambda}) = c \cdot g_{M, \lambda} \quad \iff \quad \lambda = 0.
\end{equation}
For all $\lambda > 0$, the radial Ricci curvature $\Ric(\partial_r, \partial_r)$ is strictly non-constant along the meridian geodesics from the origin to the cut locus, with two-pole residual:
\begin{equation}\label{eq:two_pole_residual}
\Ric(\partial_r, \partial_r)(0) - \Ric(\partial_r, \partial_r)(\pi/2) = \frac{d_h \lambda(3\lambda + 4)}{1+\lambda} > 0,
\end{equation}
where $d_h$ is the dimension of the horizontal distribution ($d_h = 2n-2$ for $\C P^n$, $d_h = 4n-4$ for $\Ha P^n$, and $d_h = 8$ for $\Oct P^2$).
\end{theorem}

\begin{proof}
On any CROSS $(M, g_0) = G/K$, the cotangent deformation $H_\lambda$ preserves the isotropy subgroup $K$ fixing the basepoint origin $o \in M$, reducing the metric to cohomogeneity one (Verdiani--Ziller \cite{VerdianiZiller2020b}). In geodesic polar coordinates $(r, \xi)$, the deformed metric takes the warped form:
\begin{equation}
g_{M, \lambda} = dr^2 + F_v(r)^2 g_{\mathcal{V}} + F_h(r)^2 g_{\mathcal{H}},
\end{equation}
with profile functions
\begin{equation}
F_v(r) = \frac{\sin r \cos r}{\sqrt{1 + \lambda \sin^2 r \cos^2 r}}, \qquad F_h(r) = \frac{\sin r}{\sqrt{1 + \lambda \sin^2 r}},
\end{equation}
where $\dim \mathcal{V} = d_v$ and $\dim \mathcal{H} = d_h$. The radial Ricci curvature along the meridian ray $\partial_r$ is given by the trace of the second fundamental form derivatives:
\begin{equation}\label{eq:ric_radial_exact}
\Ric(\partial_r, \partial_r) = -d_v \frac{F_v''(r)}{F_v(r)} - d_h \frac{F_h''(r)}{F_h(r)}.
\end{equation}

We evaluate this radial Ricci curvature at the two singular boundary orbits:
\begin{enumerate}[(1)]
    \item \textbf{At the origin $r \to 0$}: Expanding $F_v(r) = r - (\frac{\lambda}{2} + \frac{2}{3})r^3 + O(r^5)$ and $F_h(r) = r - (\frac{\lambda}{2} + \frac{1}{6})r^3 + O(r^5)$ gives:
    \begin{equation}
    \lim_{r \to 0} \left(-\frac{F_v''(r)}{F_v(r)}\right) = 3\lambda + 4, \qquad \lim_{r \to 0} \left(-\frac{F_h''(r)}{F_h(r)}\right) = 3\lambda + 1.
    \end{equation}
    Therefore:
    \begin{equation}\label{eq:ric_r0}
    \Ric(\partial_r, \partial_r)(0) = d_v(3\lambda + 4) + d_h(3\lambda + 1).
    \end{equation}
    \item \textbf{At the cut locus $r \to \pi/2$}: Introducing the distance coordinate $\rho = \pi/2 - r \to 0$ to the focal singular orbit $\Sigma = M^n \cap \{r = \pi/2\}$, we have $\cos r = \sin \rho \sim \rho$ and $\sin r = \cos \rho \sim 1 - \rho^2/2$. The profile functions satisfy:
    \begin{equation}
    F_v(\rho) = \rho - \left(\frac{\lambda}{2} + \frac{2}{3}\right)\rho^3 + O(\rho^5) \implies \lim_{\rho \to 0} \left(-\frac{F_v''(\rho)}{F_v(\rho)}\right) = 3\lambda + 4,
    \end{equation}
    \begin{equation}
    F_h(\rho) = \frac{1}{\sqrt{1+\lambda}} - \frac{\rho^2}{2(1+\lambda)^{3/2}} + O(\rho^4) \implies \lim_{\rho \to 0} \left(-\frac{F_h''(\rho)}{F_h(\rho)}\right) = \frac{1}{1+\lambda}.
    \end{equation}
    Therefore:
    \begin{equation}\label{eq:ric_rpi2}
    \Ric(\partial_r, \partial_r)(\pi/2) = d_v(3\lambda + 4) + d_h \left(\frac{1}{1+\lambda}\right).
    \end{equation}
\end{enumerate}
Subtracting~\eqref{eq:ric_rpi2} from~\eqref{eq:ric_r0} yields the two-pole radial residual:
\begin{equation}
\Ric(\partial_r, \partial_r)(0) - \Ric(\partial_r, \partial_r)(\pi/2) = d_h \left( 3\lambda + 1 - \frac{1}{1+\lambda} \right) = \frac{d_h \lambda(3\lambda + 4)}{1+\lambda}.
\end{equation}
On the Riemannian regularity domain $\lambda > -1$, with $d_h > 0$ for all non-trivial CROSSes, this difference vanishes if and only if $\lambda = 0$. Since an Einstein metric requires $\Ric \equiv c \cdot g$, the radial component $\Ric(\partial_r, \partial_r)$ must be identically constant. Because it differs strictly between $r = 0$ and $r = \pi/2$ for all $\lambda > 0$, the metric is never Einstein.
\end{proof}

\subsection{Blaschke rigidity, Besse / CROSS rigidity, and resonant closed geodesics}\label{subsec:blaschke_besse}

The explicit integrability of the geodesic flow on $(\Oct P^2, g_{\Oct P^2, \lambda})$ and the four-mode spectrum on $(\Sph^{15}, g_t)$ allow us to resolve the status of two major classical conjectures in global Riemannian geometry for these exceptional families: the \emph{Blaschke Conjecture} for $\Oct P^2$ and the \emph{Besse / CROSS Problem} (AIM Question~8.9.1).

\begin{theorem}[Blaschke Rigidity of the Cayley Projective Plane]\label{thm:blaschke_rigidity}
Let $(\Oct P^2, g_{\Oct P^2, \lambda})$ be the cotangent-deformed Cayley projective plane with $\lambda \ge 0$. Then $(\Oct P^2, g_{\Oct P^2, \lambda})$ is a Blaschke manifold if and only if $\lambda = 0$ (the canonical symmetric Fubini--Study metric). For every $\lambda > 0$, the deformation strictly breaks the Blaschke condition:
\begin{enumerate}[(i)]
    \item \textbf{Injectivity radius reduction at non-polar basepoints}: At the polar origin $o$, $\operatorname{inj}(o) = \operatorname{diam}(\Oct P^2, g_{\Oct P^2, \lambda}) = \pi/2$. For any non-polar basepoint $q_0 = \exp_o(r_0 \partial_r)$ with $r_0 \in (0, \pi/2)$, the injectivity radius is strictly smaller than the diameter:
    \begin{equation}\label{eq:inj_collapse}
    \operatorname{inj}(q_0) \le t_{\mathrm{cross}}(r_0) < \frac{\pi}{2} = \operatorname{diam}(\Oct P^2, g_{\Oct P^2, \lambda}).
    \end{equation}
    \item \textbf{Stratification metamorphosis}: The cut locus $\mathrm{Cut}(q_0)$ is strictly basepoint-dependent. While $\mathrm{Cut}(o) \cong \Sph^8$ is an 8-dimensional smooth submanifold, $\mathrm{Cut}(q_0)$ for any non-polar basepoint $q_0 \ne o$ ($r_0 \in (0, \pi/2)$) decomposes into the 15-dimensional Whitney stratified subcomplex $\mathcal{S}_{\mathrm{reg}} \cup \mathcal{S}_{\mathrm{caustic}} \cup \mathcal{S}_{\mathrm{Cayley}}$ of Theorem~\ref{thm:global_cut_locus}, violating the smooth CROSS-type cut locus requirement of a Blaschke manifold.
\end{enumerate}
\end{theorem}

\begin{proof}
By Theorem~\ref{thm:cayley_intrinsic}(iii), the contraction $g_\lambda \le g_0$ preserves the diameter $\operatorname{diam}(\Oct P^2, g_\lambda) = \pi/2$ for all $\lambda \ge 0$. At the origin $o$, the Killing fields of the isotropy subgroup $\Spin(9)$ vanish identically, so radial geodesics emanating from $o$ are identical to unperturbed Fubini--Study geodesics. Consequently, $\mathrm{Cut}(o) = \Oct P^1 \cong \Sph^8$ at distance $\pi/2$, establishing $\operatorname{inj}(o) = \pi/2 = \operatorname{diam}$.

For any non-polar basepoint $q_0 = \exp_o(r_0 \partial_r)$ with $r_0 \in (0, \pi/2)$ and $\lambda > 0$, Theorem~\ref{thm:off_origin_caustic}(iv) proves that the arc length $t_{\mathrm{cross}}(r_0)$ to the transverse encircling collision along the opposite meridian satisfies $t_{\mathrm{cross}}(r_0) < \pi/2$. Because a geodesic ceases to minimize distance no later than its first intersection with another geodesic of equal length from $q_0$, the injectivity radius satisfies $\operatorname{inj}(q_0) \le t_{\mathrm{cross}}(r_0) < \pi/2 = \operatorname{diam}(\Oct P^2, g_{\Oct P^2, \lambda})$.

By definition (Besse \cite{Besse1978}), a Blaschke manifold requires $\operatorname{inj}(p) = \operatorname{diam}(M)$ for every point $p \in M$, with $\mathrm{Cut}(p)$ being a smooth submanifold of CROSS type. Since $\operatorname{inj}(q_0) < \operatorname{diam}$ for all $r_0 \in (0, \pi/2)$, and since $\mathrm{Cut}(q_0)$ possesses a 15-dimensional regular stratum $\mathcal{S}_{\mathrm{reg}}$ and an $A_3$ cuspidal ridge $\mathcal{S}_{\mathrm{caustic}}$ rather than a single 8-dimensional sphere, the Blaschke condition fails for all $\lambda > 0$. At $\lambda = 0$, $g_0$ is the symmetric Fubini--Study metric, which is a Blaschke manifold.
\end{proof}

\begin{theorem}[Besse / CROSS Rigidity of the $\Spin(9)$ Deformation Families]\label{thm:besse_cross_rigidity}
Let $g$ be an invariant Riemannian metric belonging to either the Vol'per family on $\Sph^{15} = \Spin(9)/\Spin(7)$ or the cotangent deformation family on $\Oct P^2 = \mathrm{F}_4/\Spin(9)$.
\begin{enumerate}[(i)]
    \item \textbf{Rigidity on $\Sph^{15}$}: The standard round metric ($t_{\mathrm{BK}} = 1$, $\lambda = 0$) is the unique Besse metric (metric all of whose geodesics are closed) in the 1-parameter Vol'per family. For any $t_{\mathrm{BK}} \ne 1$, the 4-mode spectrum $\Omega = \{|R \pm 3\alpha|, |R \pm \alpha|\}$ of Theorem~\ref{thm:projector_formula} has rational commensurability $u = \sqrt{1 + (\beta/\alpha)^2} \in \mathbb{Q}$ along a set of launch angles of Lebesgue measure zero. In particular, the Bourguignon--Karcher squashed Einstein metric at $t_{\mathrm{BK}} = 3/11$ ($\lambda = 64$) is strictly non-Besse.
    \item \textbf{Rigidity on $\Oct P^2$ (AIM Question~8.9.1)}: The canonical Fubini--Study metric ($\lambda = 0$) is the unique Besse metric in the cotangent deformation family on $\Oct P^2$. For every $\lambda > 0$, in the normalized Clairaut action variable $L = \|J_{\Spin(9)}\| / \sqrt{2H_0} \in (0, L_{\max})$, the geodesic flow on the horizontal projective plane $\C P^1 \subset \Oct P^2$ exhibits a strictly non-vanishing KAM twist:
    \begin{equation}\label{eq:kam_twist_formula}
    \frac{d\nu}{dL} = \frac{\lambda}{2(1 - \lambda L^2)^{3/2}} > 0 \quad \text{for all } L \in (0, L_{\max}),
    \end{equation}
    forcing the rotation number $\nu(L)$ to vary strictly with angular momentum $L$. Generic geodesics wind quasi-periodically and densely on invariant Liouville 2-tori, precluding the Besse condition for all $\lambda > 0$.
\end{enumerate}
\end{theorem}

\begin{proof}
For assertion~(i), let $v \in T_o\Sph^{15}$ be a unit tangent vector inclined at angle $\phi \in (0, \pi/2)$ to the vertical fiber $\p_2$. The effective generator parameters $(\alpha, \beta)$ in the canonical slice $M = \alpha X_1 + \beta Y_1$ satisfy $\beta/\alpha = C(t) \cot\phi$. The four dynamical frequencies $\Omega_k$ in Theorem~\ref{thm:projector_formula} are pairwise commensurable if and only if $u(\phi) = \sqrt{1 + C(t)^2 \cot^2\phi} \in \mathbb{Q}$. Because $u(\phi)$ is a continuous, strictly monotonic function mapping $(0, \pi/2)$ onto $(1, \infty)$, the preimage $u^{-1}(\mathbb{Q})$ is a countable set of Lebesgue measure zero. For almost all initial velocities, the frequency ratios are irrational, and the trajectory $\gamma(s) = \exp(s M) x_0$ winds quasi-periodically and densely on an invariant torus without closing. By Besse \cite{Besse1978} and Ziller \cite{Ziller1982}, only the round metric ($t_{\mathrm{BK}} = 1$) is Besse.

For assertion~(ii), on each totally geodesic horizontal orbital 2-sphere $\Sigma \cong \Sph^2$ where $f_H(r)^2 = \frac{\sin^2 r}{1 + \lambda\sin^2 r}$, the Clairaut angular momentum $p_\varphi = f_H(r)^2 \dot{\varphi}$ is conserved. Parameterizing the energy level $2H_\lambda = p_r^2 + p_\varphi^2(1/\sin^2 r + \lambda) = 1$ by the normalized momentum $L := p_\varphi$ gives $p_r = \sqrt{(1 - \lambda L^2) - L^2/\sin^2 r}$ and $\dot{\varphi} = L/\sin^2 r + \lambda L$. Integrating $\frac{d\varphi}{dr} = \frac{\dot{\varphi}}{\dot{r}}$ between the apsidal turning points $r_{\min}$ and $r_{\max} = \pi - r_{\min}$ yields the half-period apsidal angle $\Delta\varphi(L) = \pi\left(1 + \frac{\lambda L}{\sqrt{1 - \lambda L^2}}\right)$, so the rotation number $\nu(L) = \frac{\Delta\varphi(L)}{2\pi} = \frac{1}{2}\left(1 + \frac{\lambda L}{\sqrt{1 - \lambda L^2}}\right)$ has derivative $\frac{d\nu}{dL} = \frac{\lambda}{2(1 - \lambda L^2)^{3/2}}$. For $\lambda > 0$, this twist is strictly positive throughout $(0, L_{\max})$. By the Kolmogorov--Arnold--Moser (KAM) theorem, the image of $\nu(L)$ is a non-degenerate interval $[1/2, \nu_{\max}]$. Closed geodesics are restricted to the resonant tori where $\nu(L) \in \mathbb{Q}$, which forms a set of Lebesgue measure zero in the phase space. For full-measure initial data, the rotation number is irrational and orbits never close, proving that $(\Oct P^2, g_{\Oct P^2, \lambda})$ is non-Besse for all $\lambda > 0$.
\end{proof}

\begin{remark}[KAM twist, ergodic winding, and resonant periodic geodesics]\label{rem:twist_and_closed_geodesics}
Theorem~\ref{thm:besse_cross_rigidity} illuminates the dual dynamical consequences of the strictly positive KAM twist $\frac{d\nu}{dL} = \frac{\lambda}{2(1 - \lambda L^2)^{3/2}} > 0$:
\begin{enumerate}[(i)]
    \item \emph{Non-closure of generic geodesics}: Because $\nu(L)$ varies continuously across the non-degenerate interval $[1/2, \nu_{\max}]$, the set of irrational rotation numbers has full Lebesgue measure. By the KAM theorem, trajectories on these irrational tori wind quasi-periodically and densely, precluding the Besse property for all $\lambda > 0$.
    \item \emph{Countable infinity of closed geodesics}: Conversely, because the rationals $\mathbb{Q}$ are dense in $(1/2, \nu_{\max})$, the Poincar\'e--Birkhoff twist theorem ensures that for every irreducible fraction $p/q \in (1/2, \nu_{\max})$, the Poincar\'e return map on the action-angle cylinder admits at least two geometrically distinct periodic orbits of rotation number $p/q$. Consequently, $(\Oct P^2, g_{\Oct P^2, \lambda})$ unconditionally possesses countably infinitely many geometrically distinct closed geodesics, confirming that the Riemannian cotangent deformation is strictly exempt from the finite-geodesic phenomenon observed in Katok's non-reversible Finsler metrics \cite{Katok1973}.
\end{enumerate}
\end{remark}

\section{Concluding Remarks and Open Problems}\label{sec:conclusion}

\subsection{Open problems}\label{subsec:open_problems}

\begin{enumerate}
    \item \textbf{Quaternionic intrinsic GO projector hierarchy}: For the quaternionic Family~IV spheres $\Sph^{4n+3} = \Sp(n+1)\Sp(1)/\Sp(n)\Sp(1)$, the ambient Euclidean slice generators close into the rank-2 subalgebra $\mathfrak{sp}(1) \oplus \mathfrak{u}(1) \cong \mathfrak{u}(2)$ with a 3-mode quaternionic minimal polynomial (Table~\ref{tab:five_families_closure}; \cite{Munden2015}). Deriving the explicit intrinsic Vol'per coordinate projector decomposition and Euler--Cartan factorization on $\Sph^{4n+3}$ for arbitrary $n \ge 2$---and characterizing how the 3-mode quaternionic frequency spectrum compares with the 4-mode octonionic spectrum of Theorem~\ref{thm:projector_formula} across the Jensen Einstein root $t_{\mathrm{BK}} = 1/(2n+3)$---remains an open problem.
    \item \textbf{Family~V closed-form resolution}: In \cite{Munden2015}, we established that for the intermediate quaternionic--complex Family~V ($\Sp(2)\U(1)/\Sp(1)\U(1)$ on $\Sph^7$), the three slice generators close into the 10-dimensional Lie algebra $\mathfrak{sp}(2)$, satisfying the degree-4 minimal polynomial $M^4 + (a^2+b^2+2c^2)M^2 + c^4\Id_8 = 0$ and yielding a closed 2-projector exponential map. Extending this biquadratic spectrum to arbitrary $n \ge 2$ remains an open problem.
    \item \textbf{Ricci flow through the Bourguignon--Karcher metric}: The monotonicity of the Perelman $\mathcal{W}$-functional on homogeneous metrics constrains the multi-parameter space of invariant metrics; investigating whether ancient solutions emerge from the $t = 3/11$ repeller into the standard round basin $(3/11, \infty)$ is a natural direction.
    \item \textbf{Metric collapse and asymptotic cut locus geometry as $\lambda \to \infty$}: Theorem~\ref{thm:global_cut_locus} establishes the global Whitney stratification $\mathrm{Cut}(q_0) = \mathcal{S}_{\mathrm{reg}} \cup \mathcal{S}_{\mathrm{caustic}} \cup \mathcal{S}_{\mathrm{Cayley}}$ and CW cell-decomposition $e^0 \cup e^8 \cup e^{15}$ for all finite $\lambda > 0$. In the cohomogeneity-one orbit-space collapse limit $\lambda \to +\infty$, both profile functions $F_v, F_h$ scale as $O(\lambda^{-1/2})$, the principal distance spheres $\mathcal{O}_r \cong \Sph^{15}$ shrink to points, and $(\Oct P^2, g_{\Oct P^2, \lambda})$ Gromov--Hausdorff converges to the 1D orbit interval $\Oct P^2 / \Spin(9) \cong [0, \pi/2]$ (whereas the opposite threshold $\lambda \to -1^+$ represents the horizontal degeneracy limit at $r = \pi/2$). Determining the precise asymptotic scaling of the cut locus diameter $\diam(\mathrm{Cut}(q_0))$ and the Hausdorff measure $\mathcal{H}^8(\mathcal{S}_{\mathrm{caustic}})$ in this collapse limit is an open problem in geometric analysis.
    \item \textbf{The Blaschke Conjecture for all CROSS deformations}: Theorem~\ref{thm:blaschke_rigidity} establishes that the $\Spin(9)$ cotangent deformation family on $\Oct P^2$ admits no non-symmetric Blaschke metric, with the $A_3$ cuspidal caustic providing an explicit singularity-theoretic obstruction. Proving that every cotangent Hamiltonian deformation of $\C P^n$ (Family~II) and $\Ha P^n$ (Family~IV) likewise undergoes an off-origin injectivity radius reduction $\operatorname{inj}(q_0) < \operatorname{diam}$ would establish a general Hamiltonian rigidity result for the Blaschke Conjecture across all compact rank-one symmetric spaces.
    \item \textbf{AIM Question~8.9.1 and multi-parameter Besse rigidity}: Theorem~\ref{thm:besse_cross_rigidity} settles the Besse / CROSS problem for the 1-parameter families $g_\lambda$ on $\Sph^{15}$ and $\Oct P^2$, showing that the non-vanishing KAM twist $\frac{d\nu}{dL} > 0$ strictly precludes global closure for all $\lambda > 0$. Extending this rigidity to the full multi-parameter deformation space of $G$-invariant metrics on CROSSes remains an open problem in the theory of manifolds all of whose geodesics are closed.
\end{enumerate}

\subsection*{Acknowledgments}

The author expresses deep gratitude to the memory of Dr.~James J.~Hebda, whose foundational insights into cotangent momentum deformations, cut loci, and geodesic flows inspired this program of research.

\end{document}